\documentclass[11pt]{article}
\catcode`\@=11 \catcode`\@=12
\usepackage[utf8]{inputenc}
\usepackage[T1]{fontenc}
\usepackage{amsmath}
\usepackage{amssymb}
\usepackage{amsthm}
\usepackage{oldgerm}
\usepackage{multicol}
\usepackage{verbatim}
\usepackage{color}
\usepackage{ulem}
\usepackage{enumitem}

\usepackage[pdftex,colorlinks,urlcolor=blue,pdfstartview=FitH]{hyperref}

\makeatletter
\newcommand{\mylabel}[2]{#2\def\@currentlabel{#2}\label{#1}}
\makeatother

\newcommand{\Rm}{\mathbb{R}}

\newcommand{\mL}{\mathcal{L}}

\newcommand{\mC}{\ensuremath{\mathcal{C}}}
\newcommand{\mR}{\ensuremath{\mathcal{R}}}

\newcommand{\mF}{\ensuremath{\mathcal{F}}}

\newcommand{\mD}{\ensuremath{\mathcal{D}}}

\newcommand{\mI}{\ensuremath{\mathcal{I}}}
\newcommand{\mJ}{\ensuremath{\mathcal{J}}}

\newcommand{\mE}{\ensuremath{\mathcal{E}}}

\newcommand{\vs}{\vspace{.2cm}}

\newtheorem{lem}{Lemma}[section]

\newtheorem{thm}{Theorem}

\newtheorem{cor}[lem]{Corollary}
\newtheorem{prop}[lem]{Proposition}
\newtheorem{defn}[lem]{Definition}

\newtheorem*{rmk*}{Remark}

\def\proof {\noindent{\sc{Proof. }}}
\def\qed {\mbox{}\hfill {\small \fbox{}} \\}
\def\lto{\longrightarrow}
\def\lmto{\longmapsto}

\def\leq{\leqslant}
\def\geq{\geqslant}

\begin{document}

\begin{center}
	\begin{huge}
		{\bf Lyapounov Functions of closed Cone Fields: from Conley Theory to Time Functions.}\\
	\end{huge}
	\vs
-----
\vs
\end{center}
\begin{small}
\begin{multicols}{2}

\noindent
Patrick Bernard
\footnote{Universit\'e  Paris-Dauphine},\\
PSL Research University,\\
\'Ecole Normale Sup\'erieure,\\
DMA (UMR CNRS 8553)\\
45, rue d'Ulm\\
75230 Paris Cedex 05,
France\\
\texttt{patrick.bernard@ens.fr}\\

\noindent
Stefan Suhr
\footnotemark[1],\\
PSL Research University,\\
\'Ecole Normale Sup\'erieure,\\
DMA (UMR CNRS 8553)\\
45, rue d'Ulm\\
75230 Paris Cedex 05,
France\\
%\texttt{stefan.suhr@ens.fr}\\

{ \noindent
Current address:\\
Fakult\"at f\"ur Mathematik,\\
Ruhr-Universit\"at Bohum\\
Universit\"atsstra\ss e 150,\\
44780 Bochum,\\
Germany\\
\texttt{stefan.suhr@rub.de}\\
}

\end{multicols}
\vs
\thispagestyle{empty}
\begin{center}
-----
\end{center}

\textbf{Abstract. } 
We propose a theory ``à la Conley'' for cone fields
using a notion of relaxed orbits based on cone enlargements,
 in the spirit of space time geometry.
 We work in the  setting
 of  closed (or equivalently semi-continuous) cone fields with singularities.
 This setting contains (for questions which are parametrization independent such as the existence of Lyapounov functions)
 the case of continuous vector-fields on manifolds, of differential inclusions, of Lorentzian metrics, and of continuous  cone fields.
 We generalize to this setting the equivalence between stable causality and the existence of temporal functions.
 We also generalize the equivalence between global hyperbolicity and the existence of a steep temporal function.

\begin{center}
-----
\end{center}

\textbf{R\'esum\'e. }
On développe une théorie à la Conley pour les champs de cones,
qui utilise une notion d'orbites relaxées basée sur les élargissements de cones
dans l'esprit de la géométrie des espaces temps.
On travaille dans le contexte des champs de cones fermés (ou, ce qui est équivalent, semi-continus), avec des singularités.
Ce contexte contient (pour les questions indépendantes de la paramétrisation, comme 
l'existence de fonctions de Lyapounov) le cas des champs de vecteurs continus, 
celui des inclusions différentielles, des métriques Lorentziennes, et des champs de 
cones continus.
On généralise à ce contexte l'équivalence entre la causalité stable
et l'existence d'une fonction temporale. 
On généralise aussi l'équivalence entre l'hyperbolicité globale et l'existence
d'une fonction temporale uniforme.

\vfill
\hrule
\vs
{The research leading to these results has received funding from the European Research Council
	under the European Union's Seventh Framework Programme (FP/2007-2013) / ERC Grant
	Agreement  307062. Stefan Suhr is supported by the SFB/TRR 191 `Symplectic Structures in Geometry, Algebra and Dynamics', funded 
by the DFG.}

\end{small}

\newpage

\tableofcontents

\hspace{1cm}

Lyapounov functions play an important role in dynamical systems. Their existence is related to basic dynamical behaviors
such as stability and recurrence. The second aspect was made precise by  Conley, who showed an equivalence
between the existence of Lyapounov functions and the absence of chain recurrence. This result was extended by
Hurley, see \cite{Hu1,Hu2}, to  non compact spaces. See also \cite{Pa} for a different point of view based on 
Mather-Fathi theory.

On the other hand the causality theory of space times studies (among other things)  time functions on Lorentzian manifold, see \cite{ms1} for example.
The existence of continuous  time functions  for smooth stably causal space times was proved in  \cite{geroch} and \cite{hawking}.
The condition of stable causality of space time is analogous to the absence of chain recurrence
in Conley's theory. 
Still in the context of smooth   space times,
the equivalence between  stable causality and the existence of a smooth temporal function (a regular Lyapounov function in the terminology of the present paper) was proved in \cite{besa3}. Motivated by 
solutions to the Einstein equations with low regularity the problem has been revisited  in  \cite{chrgra}, \cite{cgm} and \cite{saem}
where continuous metrics are studied.
The existence of smooth time functions 
for continuous stably causal cone fields (hence in particular for continuous, stably causal, Lorentzian metrics)
was proved in \cite{fasi} and \cite{Fa2}  by methods inspired by weak KAM theory.

In the present paper, we propose a theory ``à la Conley'' for cone fields.
Such a program was already carried out in  \cite{Mo} in the case of Lorentzian metrics,
but our  approach is different. We use a notion of relaxed orbits based on cone enlargements,
in the spirit of space time geometry.
This notion has the advantage of not resting on the choice of an auxiliary metric and
it  bypasses some technical difficulties related to the non continuity of the length.  
It allows us to work without  difficulty in the very general setting
of  closed (or equivalently semi-continuous) cone fields with singularities.
This setting contains (for questions which are parametrization independent such as the existence of Lyapounov functions)
the case of continuous vector-fields on manifolds, of differential inclusions, of Lorentzian metrics, and of continuous cone fields.
We impose a  manifold structure on the phase space, 
and  directly deal with smooth Lyapounov (or time) functions.
We generalize to this setting the equivalence between stable causality and the existence of temporal functions.
We also  prove that every globally hyperbolic cone field admits a steep Lyapounov function
(hence a Cauchy time function). The term {\it steep temporal function} 
was introduced in \cite{musa}, see  section \ref{secglhy} for the definition and a discussion.
We finally recover classical statements on the relation between Lyapounov functions and asymptotic stability
in their most general setting, as  obtained in   \cite{CLS,ST,ST2}.
Since our original motivation was to prove the existence of steep temporal functions in a generalized setting,
we work with the usual convention of space time geometry and consider Lyapounov functions which are non decreasing
along orbits (here called causal curves).

We thank the anonymous referees whose careful reading helped us writing a much better second version of the paper. Just before sending this second version, we received the 
preprint \cite{mi17} of Minguzzi, where several of our results are recovered using 
more traditional constructions.

\section{Introduction}

We work on a complete Riemannian manifold $M$.

A  {\it convex cone} in the vector space $E$ is a convex subset $C\subset E$ such that $tx\in C$ for each $t> 0$ and $x\in C$.
The convex cone $C$ is called  {\it regular} if it is not empty and it is contained in an open half-space, or equivalently
 if there exists a linear form $p$ on
$E$ such that $p\cdot v>0$ for each  $v\in C$.
The full cone $C=E$ will be called the  {\it singular cone}.
In order to shorten expressions in the sequel, we make the following definition:

\begin{defn}
	We say that $\Omega\subset E$ is an \textrm{open cone} if it is a convex cone which is open as a subset of $E$.
	
	We say the  $C\subset E$ is a \text{closed cone} if it is a convex cone which is either singular or regular and if
	$C\cup \{0\}$ is a closed subset of $E$.
	\end{defn}

Note that the empty set is both an open and a closed (regular) cone.
The empty set will be referred to as {\it degenerate}.
Note also that regular closed cones do not contain the origin.
Given $\Omega\subset E \setminus\{0\}$, we denote by $\hat \Omega$ the smallest closed cone containing $\Omega$,
and call it the closed hull of $\Omega$. 
Our definition of closed cones does not include the case of a closed half space, 
so the closed hull of an open half space is the full space.

A cone field $\mC$ on the manifold $M$ is a subset of the tangent bundle $TM$ such that $\mC(x):=T_xM \cap \mC$
is a convex cone for each $x$.
We shall only use open and closed cone fields:

\begin{defn}\label{def1}
	We say that $\mE\subset TM$ is an {open cone field} if it is a  cone field  which is open as a subset of $TM$.
	Then $\mE(x)$ is an open cone for each $x$.
	
	We say the  $\mC\subset TM$ is a {closed cone field} if it is  a cone field such that $\mC\cup T_0M$ (the zero section)
	is a closed subset of $TM$ and such that $\mC(x)$ is a closed cone for each $x$.
\end{defn}

Given a closed  cone field $\mC$, each point $x\in M$ is  of one and only one of the following types:
\begin{itemize}
\item[$\cdot$] Regular, which means that $\mC(x)$ is a regular cone, or
\item[$\cdot$] singular, which means that $\mC(x)=T_xM$, or
\item[$\cdot$] degenerate, which means that $\mC(x)$ is empty.
\end{itemize}

The  {\it domain} of $\mC$ is the set of non degenerate points. It is denoted by $\mD(\mC)$.
The domain of a closed (or open) cone field is closed  (or open).
A cone field is  called {\it non degenerate} if all points are non degenerate, \textit{i.e.} if $\mD(\mC)=M$.
The set of singular points of a closed (or open) cone field is closed  (or open).

As a first example of a closed cone field, we can associate to each continuous 
 vector field $V$ on $M$ the closed 
cone field $\mC_V$ such that $\mC_V(x)$ is the open  half line directed by $V(x)$ if $V(x)\neq 0$ and 
$\mC_V(x)=T_xM$ if $V(x)=0$. With our definitions (and this example motivates them), the singular points
of the cone field $\mC_V$ are the same as the singular points of the vector field $V$.

It is easy to see that continuous cone fields as considered for example in \cite{fasi}
are closed, hence our setting is more general.
In particular, the cone field of future directed causal vectors associated to a time oriented continuous Lorentzian
metric is a closed cone field.

The standard example of open cone field is the cone field of future directed timelike vectors 
associated to a time  oriented Lorentzian metric.

Given an immersion  $\phi: N\lto M$ and a closed (or open) cone field $\mC$ on $M$,
the pull back $\phi^*\mC:= (T\phi)^{-1}(\mC)$ is a closed (or open) cone field.
Note that the pull back may contain degenerate points even if  $\mC$ does not. This is one motivation 
to allow degenerate points.

We say that the cone field $\mC'$ is {\it wider} than the cone field $\mC$
if $\mC\subset \mC'$.
  We say that $\mC'$ is an  enlargement of $\mC$ (written $\mC\prec\mC'$)
  if there exists an open cone field $\mE$ and a closed cone field $\mD$ such that 
  $\mC\subset  \mD\subset  \mE\subset  \mC'$. An open enlargement of a closed cone field $\mC$ is just 
  an open cone field wider than $\mC$.
  
  Note that the intersection of a family of closed cone fields is a closed cone field.

  \begin{defn}
  	 Let $\mE$ be a cone field. We denote by $\hat \mE$ the smallest closed cone field  containing $\mE$, we call it the closed hull of $\mE$.
  	\end{defn}

%{\color{blue}It is easy to see that there exist nowhere singular open cone fields $\mE$ such that $\hat\mE$ is everywhere singular. The set theoretic closure 
%$\overline\mE$ in $TM$ of an open cone field $\mE$ is a cone field, though in general not a closed cone field. The difference lies in the closed set of points 
%$p\in M$ where $\overline\mE(p)$ is a half space. The exclusion of half spaces as closed cones is dictated by the generality of the theory developed here. Including closed half 
%spaces would require more restrictive assumptions. The case of cone fields with half spaces can be incorporated into the present situation as follows: Assume that $\mD$ is 
%a cone field that is closed as a subset of $TM$. Then the set of points $p\in M$ for which $\mD(p)$ is a half space is closed. Modify $\mD$ to a cone field $\hat\mD$ that is
%singular at these points $p$. Then $\hat\mD$ is a closed cone field in the sense of Definition \ref{def1}.} 
 
Given an open cone field $\mE$, we say that the curve $\gamma:I\lto M$ is {\it $\mE$-timelike}  (or just {\it timelike}) if
it is piecewise smooth (we shall see later that this regularity can be relaxed) and  if $\dot \gamma(t) \in \mE(\gamma(t))$ 
for all $t$ in $I$. At non smooth points, the inclusion is required to hold for left and right differentials.
The  {\it chronological future}  $\mI_{\mE}^+(x)$ of $x$ is the set of points $y\in M$ such that there exists a non constant  timelike curve $\gamma:[0,T]\lto M$
satisfying $\gamma(0)=x$ and $\gamma(T)=y$.
The {{\it chronological past} $\mI^-_{\mE}(x)$ of $x$} is the set of points $x'\in M$ such that $x\in \mI^{+}_{\mE}(x')$.
Note that $\mI^-_{\mE}(x)=\mI^{+}_{-\mE}(x)$.
More generally, for each subset $A\subset M$, we denote by $\mI^{\pm}_{\mE}(A):= \cup_{x\in A}\mI^{\pm}_{\mE} (x)$ the {  {\it chronological future and past
of}} $A$. They are open subsets of $M$ {  by Lemma \ref{localcone}}.
We have the inclusion $\mI_{\mE}^+(y)\subset \mI_{\mE}^+(x)$ if 
$y \in \mI_{\mE}^+(x)$.

Given a closed cone field $\mC$, we say that the curve  $\gamma:I\lto M$ is {\it $\mC$-causal}  (or just {\it causal})
if
it is locally Lipschitz and   if the inclusion $\dot \gamma(t) \in \mC(\gamma(t))\cup  
T_0 M$ 
holds for almost all $t\in I$.
The  {\it causal future}  $\mJ_{\mC}^+(x)$ of $x$ is the set of points $y\in M$ such that there exists a (possibly constant) causal curve  $\gamma:[0,T]\lto M$
satisfying $\gamma(0)=x$ and $\gamma(T)=y$.
The {  {\it causal past} $\mJ_{\mC}^-(x)$ of $x$} is the set of points $x'\in M$ such that $x\in \mJ^{+}_{\mC}(x')$.
More generally, for each subset $A\subset M$, we denote by $\mJ^{\pm}_{\mC}(A):= \cup_{x\in A}\mJ^{\pm}_{\mC} (x)$ the {\it causal future and past of $A$}. 
We have the inclusion $\mJ_{\mC}^+(y)\subset \mJ_{\mC}^+(x)$ if 
$y \in \mJ_{\mC}^+(x)$.

\begin{defn}\label{deflyap}
Let $\mC$ be a cone field on $M$. The  function $\tau:M\lto \Rm$ is called a 
Lyapounov function for the  cone field $\mC$ if it is smooth, 
$d\tau_x\cdot v\geq 0$ for each $(x,v)\in \mC$, and if,
at each regular point $x$ of $\tau$ (i.e. $d\tau_x\neq 0$),
we have $d\tau_x\cdot v>0$ for each $v\in \mC(x)$.
\end{defn}

It could be useful (especially with an eye towards degenerations of Lorentzian metrics) to study Lyapounov functions for cone fields which are closed as subsets of 
$TM$ and contain half spaces.
To a certain extent, this case can be done as follows :  Assume that $\mD$ is 
a cone field that is closed as a subset of $TM$. Then the set $\Sigma$ of points $x\in M$ for which $\mD(x)$ is a half space is closed. We can
 modify $\mD$ to the  cone field $\mC$ that is singular on $\Sigma$ and equal to $\mD$
 outside of $\Sigma$. This is a closed cone field in the sense of Definition \ref{def1}
, and the Lyapounov functions for $\mD$ are the same as the Lyapounov functions
 for $\mC$.

If $\tau$ is a Lyapounov function for the closed cone field $\mC$ on $M$, and if $\phi:N\lto M$
is an immersion, than $\tau \circ \phi$ is a Lyapounov function for $\phi^*\mC$ on $N$.

When $\mC$ is the cone field associated to a vector field $V$, a Lyapounov function for $\mC$ is the same as a Lyapounov function for $V$.

Note that if the cone field is induced by a time orientable Lorentzian metric a  Lyapounov function without critical points is a temporal function for 
the Lorentzian metric. In the same vein time/temporal function were considered in \cite{fasi} for continuous cone fields. 

Given a closed cone field $\mC$, we define
$$
\mF_{\mC}^+(x):= \{x\}\cup \bigcap_{\mE\succ\mC } \mI_{\mE}^+(x)
$$
where the intersection is taken on all open enlargements $\mE$ of $\mC$. We call $\mF_{\mC}^+(x)$ the {\it stable future} of $x$.
A point $x$ is said to be {\it stably recurrent} (for $\mC$) if, for each open 
enlargement $\mE$ of $\mC$, there exists a closed $\mE$-timelike curve
passing through $x$.
 We denote by $\mR_{\mC}$ the set of stably recurrent points.
Let us state our first result, which will be proved in Section \ref{sec-ex}.

\begin{thm} \label{thm1}
Let $\mC$ be a closed cone field.
\begin{itemize}
\item[(a)] The set $\mF_{\mC}^+(x)$ is the set of points $x'\in M$ such that $\tau(x')\geq \tau(x)$ for each (smooth) Lyapounov function $\tau$ (it is thus a closed 
set).
\item[(b)] The point $x$ is stably recurrent if and only if all (smooth) Lyapounov functions $\tau$ satisfy $d\tau_x=0$ (hence $\mR_{\mC}$ is closed).
\end{itemize}
\end{thm}

Two points $x$ and $x'$ of $\mR_{\mC}$ are called  {\it stably equivalent} if $x'\in \mF^+_{\mC}(x)$ and $x\in \mF^+_{\mC}(x')$.
This is an equivalence relation on $\mR_{\mC}$. The classes of this equivalence relation are called  {\it stable classes}.
The following statement  is also proved in  Section \ref{sec-ex}.

\begin{thm}\label{thm2}
	Let $\mC$ be a closed cone field.
There exists a (smooth) Lyapounov function $\tau$ with the following properties:
\begin{itemize}
\item[(a)] The function $\tau$ is regular at each point of $\mD(\mC)-\mR_{\mC}$.
\item[(b)] Two points $x$ and $x'$ of $\mR_{\mC}$ belong to the same stable class if and only if $\tau(x')=\tau(x)$.
\item[(c)] If $x$ and $x'$ are two points of $M$ such that $x' \in \mF^+_{\mC}(x)$ and $x \not\in \mF^+_{\mC}(x')$, then $\tau(x')>\tau(x)$.
\end{itemize}
This implies that $\mR_{\mC}$ is a closed set, as well as the stable components.
\end{thm}

We recover
the classical fact that  a closed cone field is stably causal 
(meaning that there is  no stably recurrent point ) if and only if 
it admits a smooth temporal function (in our terminology, a  Lyapounov function without critical points).
This result has a long history and several variants, see \cite{Sa05}
for the state of the art in 2005.
To our knowledge, the most general known variant before the present paper is due to
 Fathi and Siconolfi in \cite{fasi}, in the context of continuous cone fields
 (in this paper,
  the function that we call smooth temporal functions are called smooth time functions).
Our statement is more general, since we allow closed (equivalently : semi-continuous) cone fields with singularities. Our proof is entirely different.

Le us finish with a description of the stably recurrent set in terms of the 
relation $\mF^+_{\mC}$ (the analogous characterization in the Lorenzian case is given in 
\cite{Sei,HS,mi08} ):

\begin{prop}
 Let $\mC$ be a closed cone field.
	The point $x\in M$ is stably recurrent if and only if $x$ is singular or there exists a point $y\neq x$ such that $y\in \mF^+_{\mC} (x)\cap \mF^-_{\mC}(x)$.
%	\begin{itemize}
%	\item $x$ is singular 
%	\end{itemize}
%	or 	
%	\begin{itemize}
%	\item there exists a point $y\neq x$ such that $y\in \mF^+_{\mC} (x)\cap \mF^-_{\mC}(x)$.
%	\end{itemize}
\end{prop}

\proof
Let $x$ be a stably recurrent point which is not singular.
 Since the set of singular points is closed, we can choose a compact neighborhood
$K$ of $x$ which has the property that, 
for some open enlargement 
$\mE_0$ of $\mC$, all $\mE_0$-timelike loops contained in $K$ are constant.
Let $B$ be the boundary of $K$.

It follows from  Theorem  \ref{thm1} that $\mF_{\mC}^{\pm}(x)=\cap_{\mE\succ \mC}\overline{\mI^{\pm}_{\mE}(x)}$.
Indeed, if $y$ does not belong to $\mF_{\mC}^+(x)$, then there  exists  a  Lyapounov function $\tau$
such that $\tau(y)<\tau(x)$. Then, the open enlargement 
$\mE:= \{d\tau >0\}$ satisfies  $y\not \in \overline{\mI^+_{\mE}(x)}\subset \{\tau\geq \tau(x)\}$.

For each enlargement $\mE$ of $\mC$ contained in $\mE_0$, 
the closed set $\overline{\mI^+_{\mE} (x)}\cap \overline{\mI^-_{\mE}(x)}$ is not empty, 
 connected, and intersects $K$, but is not contained in $K$. 
Hence it intersects the compact set $B$.
We consider the family of non empty  closed sets 
$F(\mE):= B\cap \overline{\mI^+_{\mE} (x)}\cap \overline{\mI^-_{\mE}(x)}$
of $B$, parameterized by open enlargements $\mE$ of $\mC$.
This family  has the finite intersection property:
any finite intersection of these sets is non-empty.
By compactness of $B$, we deduce that the intersection 
$\cap_{\mE\succ \mC}F(\mE)$
on all open enlargements of $\mC$ is not empty.
This implies that $B \cap \mF^+_{\mC} (x)\cap \mF^-_{\mC}(x)$ is not empty, hence that 
$\mF^+_{\mC} (x)\cap \mF^-_{\mC}(x)$ contains a point different from $x$.
\qed

We now present some more specific applications of our methods:

\subsection{Hyperbolic cone fields}\label{secglhy}

Following the  terminology of space time geometry,
we say that the closed cone field  $\mathcal{C}$ on  $M$ is  {\it globally hyperbolic} if
\begin{itemize}
\item[\mylabel{GH0}{(GH0)}] $\mC$ is non degenerate.
\item[\mylabel{GH1}{(GH1)}] $\mathcal{C}$ is causal, i.e. all closed Lipschitz $\mC$-causal curves are constant, and that $\mC$ does not have singular points.
\item[\mylabel{GH2}{(GH2)}] The set $\mJ_{\mC}(K,K'):= \mathcal{J}_{\mC}^+(K)\cap \mathcal{J}_{\mC}^-(K')$  is compact
for each compact sets $K$ and $K'$.
\end{itemize}

We say that the closed cone field $\mC$ is {\it hyperbolic} if it satisfies \ref{GH1} and \ref{GH2}.
We stress that neither stable causality nor strong causality is assumed here,
as it is e.g. in \cite{fasi}
(it will be indirectly proved to be a consequence of hyperbolicity).
In the classical context of Lorentzian metrics,
the definition was given in a weaker form where
\ref{GH2} is replaced by

\begin{itemize}
	\item [\mylabel{GH3}{(GH3)}] The set 
	$\mJ_{\mC}(x,y)=\mJ^+_{\mC}(x) \cap \mJ_{\mC}^-(y)$ is compact for each 
	$x$ and $y$ in $M$.
\end{itemize}

Our definition is equivalent in the Lorentzian case, as follows from:

\begin{prop}\label{hypclosed}
If the closed cone field $\mC$ is wider than a non degenerate open cone field and satisfies 
\ref{GH3}, then it satisfies \ref{GH2}.
\end{prop}

\proof
Our assumption is that there exist a non degenerate open cone field $\mE\subset  \mC$. It follows from Lemma \ref{lemvf} below that $\mE$ contains a smooth vector field 
$V(x)$. This vector field can be assumed complete by reparameterization. We denote by $\phi^t$ its flow.

Let $K$ and $K'$ be two compact sets. We consider a sequence $z_n\in \mJ_{\mC}(K,K')$, 
i.e. there exist $x_n \in K$ and $y_n \in K'$ with $z_n\in \mJ_{\mC}(x_n,y_n)$. We can assume that the sequences $x_n$ and $y_n$ have limits
$x$ and $y$ in $K$ and $K'$, respectively.
For each $t>0$, 
 $x\in \mI^+_{\mE}(\phi^{-t}(x))\subset \mJ^+_{\mC}(\phi^{-t}(x))$ and 
$y \in \mI^-_{\mE}(\phi^t(y))\subset \mJ^-_{\mC}(\phi^t(y))$.
Since $\mI^+_{\mE}(\phi^{-t}(x))$ and $\mI^-_{\mC}(\phi^t(y)) $ are open,
$x_n \in \mI^+_{\mE}(\phi^{-t}(x))\subset \mJ^+_{\mC}(\phi^{-t}(x))$ and 
$y_n \in \mI^-_{\mE}(\phi^t(y))\subset \mJ^-_{\mC}(\phi^t(y))$ when $n$ is large enough,
hence $z_n \in \mJ_{\mC}(\phi^{-t}(x),\phi^t(y))$, which is a compact set
by \ref{GH3}.
We can thus assume by taking a subsequence that $z_n$ has a limit 
$z$ which is contained in  $\mJ_{\mC}(\phi^{-t}(x),\phi^t(y))$
for each $t>0$.
By \ref{GH3}, the set $\mJ_{\mC}(\phi^{-1}(x),z)$ is compact 
and it contains $\phi^{-t}(x)$ for each $t\in ]0,1[$, hence it contains $x$.
This implies that $z\in \mJ^+_{\mC}(x)$.
We prove similarly that $z\in \mJ_{\mC}^-(y)$.
\qed

The Lyapounov function $\tau$ is said to be {\it steep} if the inequality
$$
d\tau_x\cdot v\geq |v|_x
$$
holds for each $(x,v)\in \mC$. 
If $\tau$ is a steep temporal function for the closed cone field $\mC$ then 
\ref{GH1} obviously holds, and moreover  $\tau \circ \gamma(I)=\Rm$ for each complete causal curve $\gamma:I\lmto M$ (see Definition \ref{compdef}).
Indeed, such a curve has infinite length in both forward and backward direction,
and the steepness of $\tau$ then implies that $\lim _{t\lto \inf I} \tau\circ \gamma =-\infty$
and $\lim _{t\lto \sup I} \tau\circ \gamma =+\infty$.

If the cone field is non degenerate, then all inextendible causal curves
are complete, see Corollary \ref{cor-complete}.
If $\tau$ is a  steep Lyapounov  function, then 
$\tau \circ \gamma(I)=\Rm$
for each inextendible curve: steep Lyapounov functions are Cauchy time functions.
All their level sets  are {\it Cauchy hypersurfaces} in the sense that every inextendible 
causal curve intersects them exactly once.

The following statement extends a classical result (see 
 \cite{musa}, \cite{mi15}) to our more general setting. It is proved in Section \ref{sec-hyp-proof} where some other characterizations 
 of global hyperbolicity are also given.

\begin{thm}\label{thm3}
The closed cone field $\mC$ is  hyperbolic if and only if 
it admits a (smooth) Lyapounov function which is steep with respect to a complete Riemannian metric. Then, the relations $\mJ_{\mC}$ and 
$\mF_{\mC}$ are identical.
\end{thm}

As a consequence, each globally hyperbolic cone field has a Cauchy temporal function.

Note that the definition of  hyperbolicity does not involve the metric.
We deduce that, if $\mC$ is  hyperbolic and if $\tilde g$ is a 
(not necessarily complete) metric, then there exists a 
 Lyapounov function which is steep with respect to $\tilde g$. 
This follows from  the theorem applied  to the complete metric $g+\tilde g$ (where $g$
is a complete metric on $M$). 
However, a temporal function which is steep with respect to a non complete metric
is not necessarily a Cauchy time function, even in the absence of degenerate points. 
The existence of such a function does not necessarily
imply hyperbolicity. As an example, consider a complete Riemannian manifold $(M,g)$, a globally hyperbolic
cone field $\mC$, and a steep temporal function $\tau$. If $N$ is the complement of a point in $M$,
then $(N,g)$ is not complete, $\tau$ is a steep time function which is not Cauchy, and $\mC$ is not  hyperbolic.

Our notion of steep Lyapounov functions is similar to, but different from, the notion of steep temporal function introduced
in \cite{musa}, see also \cite{mi17}, as the sharp criterion for the isometric 
embeddability of space times into Minkowski space. 
 There a function $\tau$ on the 
space time $(M,g_L)$ is called  steep if $d\tau\cdot v\ge \sqrt{|g_L(v,v)|}$ for all 
future pointing vectors $(x,v)\in TM$. 
Since we can choose a metric $g$ such that  $g(v,v)\ge |g_L(v,v)|$ on all tangent vectors,
 Theorem \ref{thm3} 
implies  the existence of a  steep temporal functions in the sense of  \cite{musa} in globally hyperbolic space times.
Conversely the existence of a steep temporal functions in the sense of \cite{musa} does not  imply global hyperbolicity.

The conclusion of Theorem \ref{thm3} is false if \ref{GH2} is  replaced by \ref{GH3} without assuming 
that $\mC$ has non empty interior.
Any vector field admitting non trivial recurrence provides a counter-example.
We have the following corollaries:

\begin{cor}
	Each  hyperbolic cone field admits a  hyperbolic enlargement. Especially globally hyperbolic cone fields have globally hyperbolic enlargements.
\end{cor}

\proof
Let $\tau$ be a steep Lyapounov function.
The closed cone field $\mC$ is contained in the closed cone field $\{(x,v): x\in \mD(\mC) \text{ and } d\tau_x\cdot v \geq |v|_x\}$.
Let $F$ be a closed set containing $\mD(\mC)$ in its interior and disjoint from the critical set of $\tau$.
The closed cone field  
$$\mathcal{G}:=\{(x,v): x\in F \text{ and } d\tau_x\cdot v \geq |v|_x/2\}$$ 
is thus an enlargement of $\mC$, and $2\tau$ is a steep temporal function for it, hence it is  hyperbolic.
If $\mC$ is globally hyperbolic $\mD(\mC)=M$ and therefore $F=M$. Then $\mathcal{G}$
is globally hyperbolic.
\qed

In particular, hyperbolicity implies stable causality and strong causality.  The definition of strong causality is analogous to 
the one in Lorentzian geometry. More precisely, a cone field $\mC$ is {\it strongly causal} if every point has a neighborhood $U$ such 
that the set $\gamma^{-1}(U)$ is connected set in the interval $I$ for every causal curve $\gamma\colon I\to M$. 

Although the notion of steep temporal functions is less intrinsic  than the notion of Cauchy 
temporal function (it depends on the choice of an auxiliary metric), the above proof 
shows that it is more tractable, compare \cite{Sa13}.

The  splitting theorem, see \cite{besa2,besa3}, also holds in our setting:

\begin{cor}\label{cor3}
Let $(M,\mathcal{C})$ be globally hyperbolic. 
Then there exists a manifold $N$ and a diffeomorphism $\psi:M\lto \Rm\times N$
whose first component is a steep time function on $M$.
\end{cor}

\proof
Let $\tau$ be a steep time function.
We consider the vector field
$V(x)=\nabla \tau/|\nabla \tau |^2$, which has the property that 
$d\tau_x \cdot V(x)=1$. Note that $|d\tau_x|\geq 1$ hence 
$|\nabla \tau_x|\geq 1$ hence $|V(x)|\leq 1$.
As a consequence, the flow $\varphi^t$  of $V$ is complete.
Setting $N=\tau^{-1}(0)$, the map $\phi:(t,x)\lmto \varphi^t(x)$ is 
a diffeomorphism from $\Rm\times N$ into $M$. 
For each point $x\in M$, we have $\varphi^{-\tau(x)}(x)\in N$, hence 
$x\in \varphi^{\tau(x)}(N)$. This implies that $\phi$ is onto, and that the 
first component of the inverse $\psi$ of $\phi$ is equal to $\tau$.
\qed

As was noticed in \cite{chenem}, if $M$ is moreover assumed contractible, 
it is then  diffeomorphic  to a Euclidean space.

\subsection{A Lemma of Sullivan}

We start with the definition of complete causal curves, which are the analogs in our setting of maximal solutions 
of vector fields.

\begin{defn}\label{compdef}
	The causal curve $\gamma$ is called complete
	if it is defined on an open (possibly unbounded) interval $]a,b[$ 
	and if the two following conditions hold:
	\begin{itemize}
		\item[(a)]
		Either $\gamma|_{[s,b[}$ has infinite length 
		for each $s\in ]a,b[$ or $\lim_{t\lto b} \gamma(t)$
		is a singular point of $\mC$
		(we say that $\gamma$ is forward complete).
		\item[(b)]
		Either $\gamma|_{]a,s]}$ has infinite length
		for each $s\in ]a,b[$  or $\lim_{t\lto a} \gamma(t)$
		is a singular point of $\mC$
		(we say that $\gamma$ is backward complete).
	\end{itemize}
\end{defn}

Although the notion of complete curve is expressed in terms of an auxiliary complete metric $g$, it is not hard to see that 
it does not depend on $g$, as long as $g$ is complete.
We have:

\begin{prop}\label{prop-sul}
	Let $(M, \mC )$ be a closed cone field and let $F\subset M$ 
	be a closed set. Let $Z\subset F$ be the union of all complete causal curves contained in $F$.
	Then, there exists a Lyapounov function $\tau$ for $\mC$ on $M$ which is regular on $F-Z$.
\end{prop}

\proof
We consider the closed cone field $\mC_F$ which is equal to $\mC$ on $F$
and degenerate outside of $F$.
Each curve which is causal and complete for $\mC_F$ is causal and complete  for $\mC$.
The proposition follows from Theorem  \ref{thm2} 
and the observation that $\mR(\mC_F)\subset Z$, which follows from 
Corollary \ref{cor-R} below, applied to $\mC_F$.
\qed

In the case where $\mC$ is the cone field generated by a continuous vector field $X$, where $F$
is compact, and where $Z$ is empty,
we recover the following  famous  Lemma of Sullivan,  \cite{sul}:

\textit{If $X$ is a continuous vector field on $M$, and if $K$ is a compact set which 
	does not contain any full orbit of $X$, then there exists a Lyapounov function 
	for $X$ which is regular on $K$, \textit{i.e.}  $d\tau_x\cdot X(x)>0$
	for each $x\in K$.}

The proof of Sullivan  in \cite{sul} was based on  the Hahn-Banach Theorem,
a  more elementary proof was given in \cite{LS}. Proposition \ref{prop-sul} extends this result to the non compact case,
and also to the case where some full orbits exist.

\subsection{Asymptotic stability}

We consider a closed cone field $\mC$.
A compact set $Y\subset M$ is called { \it asymptotically stable} if, for each neighborhood
$U$ of $Y$, there exists a neighborhood $V\subset U$ of $Y$ such that 
$\mJ^+_{\mC}(V)\subset U$ and if each forward complete causal curve starting in $V$ converges to $Y$
(which means that the distance to $Y$ converges to zero).
If $Y=\{y\}$ is a point, then this requires that $y$ be singular (or degenerate).

We can recover in our setting the following restatements of several known results 
on converse Lyapounov theory for differential inclusions, see \cite{CLS} for the case where $Y$ is a singular point,
and \cite{ST,ST2} for the general case. Our setting in terms of cone fields is parametrization-invariant,
in contrast to the formulation in terms of differential inclusions used in the papers mentioned before. 
Since both  properties of being asymptotically stable and of admitting a Lyapounov function 
are parametrization invariant, these settings are equivalent.
Note that our sign convention for Lyapounov functions is non standard: They increase along orbits.

\begin{prop}
	Let $Y\subset M$ be a compact set and let $\mC$ be a closed cone field which is 
	non degenerate in a neighborhood of $Y$. The following properties are equivalent:
\begin{enumerate}
	\item \label{AS1} $Y$ is asymptotically stable.
	\item \label{AS2} $\mJ_{\mC}^+(Y)=Y$ and there exists a  neighborhood $U$ of $Y$ such that each backward complete causal curve $\gamma$
	contained in $U$ is contained in $Y$.
	\item  \label{AS3} $\mF_{\mC}^+(Y)=Y$ and there exists a compact  neighborhood $U$ of $Y$
	such that $\mR_{\mC}\cap U\subset Y$.
	\item  \label{AS4} There exists a  Lyapounov function $\tau$ which is  null on $Y$,
	regular  on $U-Y$, and negative on $U-Y$,  where $U$ is a neighborhood of $Y$.
\end{enumerate}
\end{prop}

\proof
$\ref{AS1}\Rightarrow \ref{AS2}$. 
The asymptotic stability implies that $\mJ_{\mC}^+(Y)\subset U$ for each neighborhood $U$ of $Y$, 
hence  $\mJ_{\mC}^+(Y)\subset Y$.
Let $U_0$ be a compact neighborhood of $Y$ which has the property that all forward complete curves contained in $U_0$ 
converge to $Y$. This implies in particular that $U_0-Y$ does not contain singular points.
Let us suppose that there exists a   backward complete causal curve 
$\gamma:]-T,0]\lto  U_0$ such that $\gamma(0)$ does not belong to $Y$. 
Let $U_1$ be a compact neighborhood of $Y$ which does not contain $\gamma(0)$.
There exists an open  neighborhood $V_1$ of $Y$ such that $\mJ_{\mC}^+(V_1)\subset U_1$, which implies that $\gamma$ does not enter $V_1$ on $]-T,0]$.
Since $U_0-V_1$ does not contain singular points of $\mC$, the curve $\gamma$ has infinite length,
we parametrize it by arclength, $\gamma:(-\infty, 0]\lto M$. By the Ascoli Arzela Theorem, there exists a sequence $t_n\lto -\infty$
such that the curves $t\lmto \gamma(t-t_n)$ converge, uniformly on compact intervals, to a Lipschitz  curve $\eta:\Rm\lto U_0-V_1$.
By Proposition \ref{proplim}, the curve $\eta$ is causal and forward complete.
This  implies that $\eta$ converges to $Y$, which is a contradiction since $\eta (\Rm) \subset U_0-V_1$.

$\ref{AS2}\Rightarrow \ref{AS3}$.
Let $U$ be the neighborhood with property \ref{AS2}, and $W$ be a compact neighborhood 
of $Y$ contained in $U$.
If $\mF_{\mC}^+(Y)$ was not contained in $W$, then there would exists a backward complete causal curve 
contained in $W$ but not in $Y$, by Corollary \ref{cor-back}.
This contradiction implies that  $\mF_{\mC}^+(Y)\subset W$, and, since this holds for each compact neighborhood 
 $W$ of $Y$ contained in $U$, that $\mF_{\mC}^+(Y)\subset Y$.
The part of the statement concerning $\mR_{\mC}$ follows immediately from Corollary \ref{cor-R}.

$\ref{AS3}\Rightarrow \ref{AS4}$. It  is a direct consequence of Proposition \ref{proplyapnulcomp}.

$\ref{AS4}\Rightarrow \ref{AS1}$ Let $U$ be a compact neighborhood of $Y$ such that $ \tau$ is regular and negative on $U-Y$.
For each compact neighborhood $W$ of $Y$ contained in $U$, we set  $a=\max_{\partial W}\tau$ (by compactness, $a<0$)
and $V:= \{x\in W, \tau(x)\geq a/2\}$.
We have $\mJ^+_{\mC} (V)\subset V\subset U$.
Let $\gamma:[0,T[\lto V$ be a complete causal curve parametrized by arclength.
The function $\tau \circ \gamma$  is non decreasing, hence it converges to $b\in [a/2,0]$.
We have to prove that $b=0$.
The set $V^b:= \{x\in V, \tau(x)\leq b\}$ is compact.
If $b<0$, then $\tau$ is regular on $V^b$, hence there exists $\delta>0$ such that 
$d\tau_x \cdot v\geq \delta |v|$ for each $(x,v)\in \mC, x\in V^b$.
This implies that $\tau\circ \gamma(t)\geq \tau \circ \gamma(0)+\delta  t$, hence that $T\leq (b-a)/\delta$ is finite.
The complete causal curve $\gamma$ has finite length, hence it converges to a limit $x\in V^b$ which is a singular point of $\mC$
hence a critical  point of $\tau$, a contradiction.
\qed

 \section{Preliminaries}

 \subsection{On cone fields}
 
 We state here useful results on cone fields.
 
 \begin{lem}
 If $\mC$ is a closed cone field and $\mE$ an open cone field, then
 the set of points $x\in M$ such that $\mC(x)\subset  \mE(x)$ is open.
 \end{lem}
 
 \proof
 It is the projection on $M$ of the open set $\mE-\mC$.
 \qed

 A standard partition of the unity argument implies:

\begin{lem}\label{lemvf}
 Let $\mE$ be an open cone field.
Then there exists a smooth vector field $V$ on $\mD(\mE)$  such that $V(x)\in \mE(x)$ for each $x\in \mD(\mE)$.
Moreover, given  $(x,v)\in \mE$, the vector field $V$ can be chosen such that 
 $V(x)=v$.
 In particular, there exists a smooth curve $\gamma(t): \Rm \lto M$ which is $\mE$-timelike and such that $(\gamma(0), \dot \gamma(0))=(x,v)$.
\end{lem}

Note that $V(x)\neq 0$  at each  regular point $x$ of $\mE$. This implies that a non degenerate open cone field on
a manifold must admit singular points if the Euler characteristic is not zero.

\begin{lem}\label{lem-infinitelengh}
Let $\mE$ be a non degenerate  open cone field.
Then from each point $x$ starts a forward timelike curve of infinite length, and a backward
timelike curve of infinite length.
\end{lem}

\proof
We consider a smooth vector field $V(x)$ contained in $\mE$. 
We assume that $V$ has a complete flow (this can  be achieved by multiplying $V$
by a positive smooth  function).
Let $\gamma(t):[0,\infty)\lto M$ be the forward orbit of $x$ under this flow.
Either $\gamma$ has infinite length, or it converges to a singular point of $\mE$.
In the second case  there exists a finite time $T>0$ such that $\gamma(T)$ is singular,
since the singular set is open.
We can then extend the curve 
$
\gamma_{|[0,T]}
$
by a small loop at $\gamma(T)$ contained in the singular set of $\mE$,
and obtain this way a timelike curve of infinite length starting at $x$. 
The backward case is analogous.
\qed

A smooth function $\tau$ defined near $x$ is called a {\it local Lyapounov function
at $x$} (for the closed cone field $\mC$) if $d\tau_x\neq 0$ and $d\tau_x\cdot v>0$ for each  
$v\in \mC(x)$.
This property then holds in a neighborhood of $x$.
 Local Lyapounov functions at $x$ exist if and only if 
$x$ is a not a singular point of $\mC$.
The cone $\mC(x)$ is the set of vectors $v\in T_xM$ such that 
$d\tau_x \cdot v > 0$ for each local Lyapounov function $\tau $ at $x$.
The set of vectors $v\in T_xM$ such that 
$d\tau_x \cdot v \geq 0$ for each local Lyapounov function $\tau $ at $x$ is $\mC(x)\cup \{0\}$.

Let $C$ be a closed cone and $\Omega \succ C$ be an open cone. Then there exists an open cone $\Omega'$
such that $C\subset  \Omega' \subset  \hat \Omega' \subset  \Omega$.
Given a diffeomorphism onto its image $\phi: N\lto U\subset M$
and a cone field $\mC$ on $M$,  we denote by $\phi^*\mC:=(T\phi)^{-1}(\mC)$
the preimage  of the cone field $\mC$, where $T\phi$ is the tangent map $(x,v) \lmto (\phi(x), d\phi_x\cdot v)$.
Similarly we define the forward image $\phi_*\mC:=T\phi(\mC)$ of a cone field on $N$, this is a cone field on $U=\phi(N)$.
We denote by $Q_s, s\geq 0$ the standard open cone 
$$Q_s:=\{(y,z)\in \Rm^{d-1}\times \Rm :z> s|y|\}\subset \Rm^d.
$$
In the following and later, we denote by $B^d(r)$ the Euclidean open ball of radius 
$r$ in $\Rm^d$, and also set $B^d:= B^d(1)$.
 
 \begin{lem}\label{localcone}
 Let   $\mE$ be an open  cone field  and  let $x_0$ be a point
 which is non degenerate for $\mE$ 
 and regular for the closed hull $\hat \mE$.
There exists a chart $\phi: B^{d-1}\times ]-1,1[\lto M$ at $x_0$ such that 
$$
Q_1\subset  \phi^* \mE(y,z)\subset  \phi^* \hat \mE (y,z)\subset  Q_0
$$
for each $(x,y)\in B^{d-1}\times ]-1,1[$.
\end{lem}

\proof
Let $\tau$ be a local  Lyapounov function for $\hat \mE$ such that $\tau(x_0)=0$.
Let $V$ be a vector contained in $\mE(x_0)$ and $\psi:M\lto \Rm^{d-1}$
be a smooth local map sending $x_0$  to $0$ and such that 
the kernel of $d\psi_{x_0}$ is $\Rm V$. For each $a>0$, the map 
$\Psi:=(a\tau, \psi)$ is  a local diffeomorphism, such that 
$d \Psi_{x_0}\cdot V=(ad\tau_{x_0}\cdot V, 0)$ and $\Psi_* \hat \mE(0,0)\subset  Q_0$. 
If $a>0$ is small enough, we have 
 $Q_{1/2} \subset  \Psi_* \hat \mE(0,0)$.
 As a consequence, there exists $s>0$ such that 
$$
Q_1\subset  \Psi_* \mE\subset  \Psi_* \hat \mE \subset  Q_0
$$
on $B^{d-1}(s)\times ]-s,s[$. The inverse map $\phi$ of $\Psi/s$ 
then satisfies the conclusions of the statement.
\qed

The following classical observation will be useful.

\begin{lem}\label{closure}
Let $\mE$ be an open cone field. For each subset $A$ of $M$, we have 
$$
\mI^{\pm}_{\mE}(\bar A)=\mI^{\pm}_{\mE}( A)
$$
\end{lem}

\proof
Let us consider a point $y\not \in \mI^+_{\mE}(A)$. Then the open set $\mI^-_{\mE}(y)$ 
is disjoint from $A$, hence from $\bar A$. We conclude that $y\not \in \mI^+_{\mE}(\bar A)$.
\qed

In the sequel we will need 
 the notion of sums of convex cones or cone fields.
The sum of a family of convex cones is defined
as the convex envelop of their union.
The sum of  cone fields is defined pointwise.

\begin{lem}
The sum $\mE=\sum _{\alpha} \mE_{\alpha}$ of an arbitrary family of open cone fields is an open cone field.
\end{lem}

\proof
Let $(x,v)\in\mE$.
We can assume that $M=\Rm^d$ by working in a chart at $x$.
 The vector $v$ belongs to the convex closure of the union $\cup_{\alpha} \mE_{\alpha}(x)$,
hence it is a finite sum of elements of this union:
There exists a finite set $J$ of indices such that
$v=\sum_{i\in J} v_i$ with $v_i \in \mE_i(x)$.
Let $B_i\subset \mE_i(x)$ be a compact neighborhood of $v_i$ in $\Rm^d$.
For each $i\in J$, there exists a neighborhood of $x$ on which
$B_i\subset  \mE_i(y)$. As a consequence, there exists a neighborhood
$U$ of $x$ such that $B_i \subset \mE_i(y)$ for each $y\in U$ and each $i\in J$.
We conclude that $U\times (\sum_i B_i) \subset \mE$.
\qed

\begin{lem}\label{sum}
Let $\Omega$ be an open cone and let $C_i$ be finitely many closed cones
such that $C_i\subset  \Omega$. Then there exists an open cone $\Omega'$ such that 
$\sum C_i \subset  \Omega' \subset  \hat \Omega'\subset  \Omega$.
\end{lem}

\proof
In the case where $\Omega =\Rm^d$, we can take $\Omega'=\Omega$.
Otherwise we can assume that $\Omega \subset  Q_0$ (the open upper half space).
Each of the closed cones $C_i$ then satisfies $C_i\subset  Q_{s_i}$ for
some $s_i>0$.
We can take $\Omega'=Q_s$ with $s=\min s_i$.
\qed

\begin{lem}\label{intermediate}
If $\mE$ is an open enlargement of the closed cone field $\mC$, then there exists an open cone field $\mE'$
such that
$$
\mC \subset  \mE' \subset   \widehat{ \mE'}\subset  \mE.
$$
\end{lem}

 \proof
For each $x_0\in M$, there exists a chart $\phi:B^{d-1}\times ]-1,1[\lto M$ at $x_0$ and an open  cone $\Omega\subset \Rm^d$
such that $\phi^* \mC(y,z)\subset  \Omega \subset  \hat \Omega \subset  \phi^*\mE (y,z)$ for each $(y,z)\in B^{d-1}\times ]-1,1[\lto M$ .
We take a locally finite covering of $M$ by open sets $U_i$ which
are of the form $\phi_i(B^{d-1}(1/2)\times ]-1/2,1/2[)$ for such charts,
and denote by $\Omega_i$ the corresponding open cones.
We consider the open cone fields $\mE_i$
which are equal to $\phi_* \Omega_i$ on $U_i$ and which are empty outside 
of $U_i$.
The closed hull $\hat \mE_i$ is the cone field equal to $\phi_* \hat \Omega_i$ on $\bar U_i$
and empty outside of $\bar U_i$.
Then we consider the open  cone field $\mE'=\sum_i \mE_i$.
For each $x\in M$, there exists $i$ such that $x\in U_i$,
hence $\mC(x)\subset  \mE_i(x)\subset  \mE'(x)$.

Let us now prove that $\widehat{\mE'}(x)\subset  \mE(x)$ for each $x\in M$.
Let $J(x)$ be the finite set of indices such that $x$ belongs to the closure of 
$U_i$.
For each $i\in J(x)$, $ \widehat{\mE_i}(x)=\hat \Omega_i\subset  \mE(x)$.
Lemma \ref{sum} implies the existence of   a convex  open cone $\Omega'\subset T_xM$ such that 
$$\widehat {\mE_i}(x)\subset \Omega'\subset  \widehat{\Omega'} \subset    \mE(x)
$$
for each $i\in J(x)$. 
We { use} a chart at $x$ to identify the tangent spaces $T_yM$ with $T_xM$ for $y$ near $x$.
The inclusion 
$$
\widehat{ \mE_i}(y)\subset  \Omega'\subset   \widehat {\Omega'}\subset  \mE(y)
$$
holds for $y$ in  an open neighborhood $V$ of { $x$}.
Let us consider the closed cone field $\mC'$ which is equal to $\widehat {\Omega'}$ on $V$ and which
is singular on $M-V$. We have
$
\widehat{ \mE_i}\subset  \mC' 
$
 for each $i$,
which implies that 
$\mE'\subset  \mC'$, and then that $\widehat{ \mE'}\subset  \mC'$.
Since $\mC'(x)\subset  \mE(x)$, we deduce that $\widehat{ \mE'}(x)\subset  \mE(x)$.
 \qed
 
 \begin{lem}\label{decrease}
There exists a sequence $\mE_n$ of open cone fields
which is strictly decreasing to $\mC$, which means that  $\widehat{ \mE_{n+1}}\subset \mE_n$ for each $n$ and that 
$\mC=\cap \mE_n$. 
Such a sequence has the property that, for each open enlargement $\mE$
of $\mC$ and each compact set $K\subset M$, there exists $n$ such that $\mE_n(x)\subset  \mE(x)$
for each $x\in K$.
 \end{lem}
 
\proof
For each point $(x,v)\in TM-{ (\mC\cup T_0M)}$, there exists an open enlargement 
$\mE$ of $\mC$ which is disjoint from a neighborhood $U$ of $(x,v)$.
We can  cover the complement of $\mC$ in $TM$ by a sequence $U_i$
of open sets such that, for each $i$, there exists an open enlargement $\mE'_i$
of $\mC$ disjoint from $U_i$.
We define inductively the open cone field  $\mE_{n}$ as an enlargement
of $\mC$ satisfying 
$$
\widehat{\mE_n}\subset  \mE'_n \cap \mE_{n-1}.
$$
It is obvious from the construction that $\mC=\cap \mE_n$.
Finally, let $K\subset M$ be compact and $\mE$ be  an open enlargement of $\mC$.
For each $x\in K$, there exists $n_x$ such that 
$\widehat {\mE_{n_x}}(x)\subset  \mE_{n_x-1}(x)  \subset \mE(x)$.
Then the inclusion $\widehat { \mE_{n_x}}(y)  \subset \mE(y)$
holds on an open neighborhood of $x$.
{ Using Lemma \ref{intermediate} we} can cover $K$ by finitely many such open sets, hence $\mE_n(y)\subset  \mE(y)$
for each $y\in K$ when $n$ is large enough.
\qed

 \subsection{Clarke differential, causal and timelike curves}

 We will use the notion of  Clarke differential of curves and functions, see \cite{C:90} for example.
 
 The Clarke differential of a locally Lipschitz function $f:\Rm\lto \Rm$ at a given point $x$
 is the compact interval 
 $$
 \partial f(x) =\left[\liminf _{y_2\rightarrow x, y_1\rightarrow x, y_2>y_1}
 \frac{f(y_2)-f(y_1)}{y_2-y_1},
 \limsup _{y_2\rightarrow x, y_1\rightarrow x, y_2>y_1}
 \frac{f(y_2)-f(y_1)}{y_2-y_1}\right].
 $$
 The interval $\partial f(x)=[p_-,p_+]$ can be characterized in the following way:
 for $p<p_-$, the function $t\lmto f(t)-pt$ is increasing near $t=x$,
 it is decreasing for $p>p_+$, and it is not monotone in any neighborhood of $x$ for 
 $p\in ]p_-,p_+[$.

 The Clarke differential of a locally Lipschitz curve $\gamma:\Rm\lto M$ at a given time $t$ 
 is the compact convex  subset  $\partial \gamma(t)\subset T_{\gamma(t)} M$ 
  defined 
  as the convex hull  of limit points of sequences of the form 
  $(\gamma(t_n),  \gamma'(t_n))$ in $TM$, where $t_n$ is a sequence of differentiability points of $\gamma$, see { \cite[Theorem 2.5.1]{C:90}}.
 It satisfies the equality 
 $$
 df_{\gamma(t)}\cdot \partial \gamma(t) =\partial(f\circ \gamma )(t)
  $$
  for each smooth function $f$, and this characterizes $\partial \gamma(t)$.
  In other words, $v\in \partial \gamma(t)$ if and only if 
  $
  df_{\gamma(t)}\cdot v \in \partial(f\circ \gamma )(t)
  $
  for each smooth function $f$.

\begin{lem}\label{lem-causal}
	Given a closed cone field $\mC$ on $M$, the following statements are equivalent 
	for a locally Lipschitz curve $\gamma:I\lto M$:
	
	\begin{itemize}
		\item[(a)]
		$ \gamma'(t)\in \mC(\gamma(t))$ for almost every $t\in I$.
		\item[(b)]
		$\partial \gamma(t) \subset \mC(\gamma(t))$ for each $t\in I$.
		\item[(c)]
		For each $t\in I$ and each local Lyapounov function $\tau$ at $\gamma(t)$, the function
		$\tau \circ \gamma$ is non decreasing in a neighborhood of $t$.
	\end{itemize}
\end{lem}
%We call the corresponding Lipschitz curves $\mC$-causal.

\proof
Note that (b) and (c) both hold at $t$ if $\gamma(t)$ is a singular point of $\mC$.
In this case the lemma is trivial. Therefore we can assume that $\gamma(t)$ is a regular point of $\mC$.

Property (b) implies (a) since $\gamma'(s)$ exists almost everywhere, and then is contained in $\partial \gamma(s)$.

Assume property (a).
For each $t\in I$ and each local time function $\tau$ at $\gamma(t)$,
we consider a neighborhood  of $\gamma(t)$ such that $\tau$
is a regular Lyapounov function on $U$.
We have $\gamma(s)\in U$ for $s$ close to $t$.
Then, for almost every point $s$ in a neighborhood of $t$, the derivative
$(\tau\circ \gamma)'(s)=d\tau_{\gamma(s)}\cdot \gamma'(s)$ exists and  is non negative.
This implies that the Lipschitz function $\tau \circ \gamma$ is non decreasing near $t$.

If (b) does not hold at some time $t$, then there exists $w\in \partial \gamma(t)$
and a local time function $\tau$ at $\gamma(t)$ such that $d\tau_{\gamma(t)}\cdot w<0$.
This implies that $\partial (\tau \circ \gamma)(t)$ contains a negative value, hence that  
$\tau\circ \gamma$ is not non decreasing near $t$, i.e. contradicting (c).
\qed

 Given an open cone field $\mE$, we call a locally Lipschitz curve {\it $\mE$-timelike} if it satisfies the inclusion
$\partial \gamma(t) \subset \mE(\gamma(t))$ for each $t$.
In the case of open cone fields, this is stronger than requiring the inclusion $\gamma'(t)\subset\mE(\gamma(t))$ 
for almost every $t$. 
Note however that a piecewise smooth curve is timelike according to the present definition if and only if 
it satisfies the definition given in the introduction
(the Clarke differential at a non smooth point of a piecewise smooth curve is the interval
whose endpoints are the left and right derivatives at that point).
It is equivalent to define the timelike future $\mI^+_{\mE}$ using smooth, piecewise smooth, or Lipschitz timelike curves :

\begin{lem}
	Let us consider an open cone field $\mE$.
	For each Lipschitz timelike curve $\gamma\colon[0,T]\lto M$, there exists a smooth timelike curve $\tilde \gamma: \Rm\lto M$ such that 
	$\tilde \gamma(0)=\gamma(0)$ and $\tilde \gamma(T)=\gamma(T)$.
\end{lem}

\proof
Let us set $x=\gamma(0)$ and work in a local chart at $x$.
Since $\partial \gamma(0)\subset \mE(x)$, there exists a compact neighborhood $K$ of $\partial\gamma(0)$ 
which is contained in $\mE(x)$.
In view of the semi-continuity of the Clarke differential, we deduce that $\dot \gamma(t)\in K$ almost everywhere in some interval $]0, \epsilon[$.
As a consequence, the curve $\gamma(0)+t(\gamma(\epsilon)-\gamma(0))/\epsilon$ is timelike on $[0, \epsilon]$.
By using the same procedure at $T$, we can assume that $\gamma$ is smooth near the boundaries.
We can then smooth $\gamma$ on the whole interval by a standard convolution to a smooth timelike curve.
\qed

We will also use the concept of Clarke differential of a locally Lipschitz function $f: M\lto \Rm$.
The Clarke differential $\partial f(x)\subset T^*_xM$ at the point $x$ is 
defined 
as the convex hull  of limit points of sequences of the form 
$(x_n,  df(x_n))$ in $T^*M$, where $x_n$ is a sequence of differentiability
 points of $f$ converging to $x$, see {\cite[Theorem 2.5.1]{C:90}}.
It satisfies the equality 
$$
\partial f_{\gamma(t)}\cdot \gamma'(t) =\partial(f\circ \gamma )(t)
$$
for each smooth curve $\gamma$, and this characterizes $\partial f$.
In other words, $p\in \partial f(x)$ if and only if 
$
p\cdot \gamma'(t)\in \partial (f\circ \gamma)(t)
$
for each smooth curve $\gamma$ satisfying $\gamma(t)=x$. 
If $h$ is a $C^1$ function, then  $\partial h(x)=\{dh(x)\}$,
and more generally {$\partial(h+f)(x)= dh(x)+\partial f(x)$} for each locally Lipschitz function $f$.

If $\gamma$ is a locally Lipschitz curve and $f$ a Lipschitz function, 
then we have the chain rule (see {\cite[Theorem 2.3.9]{C:90}})
$$
\partial (f\circ \gamma)(t_0)\subset [\inf_{p,v}p\cdot v, \sup_{p,v} p\cdot v]
$$
where the $\sup$ and $\inf$ are taken on  
$p\in \partial f(\gamma(t_0)),
v\in \partial \gamma(t_0)$.
This inclusion is an equality if $\gamma$ or $f$ are smooth, as seen above,  but may be strict 
in general.

\subsection{Limit curve Lemma}

Recall that we have fixed a complete Riemannian metric on $M$.
We consider a closed cone field $\mC$ and  a sequence $\mE_n$ of open cone fields strictly decreasing to 
$\mC$ in the sense of Lemma \ref{decrease}.

\begin{lem}\label{lem-limit}
	Let $\gamma_n:I\lto M$ be an equi-Lipschitz sequence of $\mE_n$-timelike curves
	converging to $\gamma:I\lto M$ uniformly on compact subintervals of $I$,
	then $\gamma$ is $\mC$-causal.
	\end{lem}

\proof
Note first that $\gamma$ is Lipschitz. 
Let $t\in I$ be given, and let $\tau$ be a local Lyapounov function at $\gamma(t)$.
In view of Lemma \ref{lem-causal}, it is enough to prove 
that $\tau \circ \gamma$ is non decreasing near $t$.

Let $U$ be a compact neighborhood of $\gamma(t)$ such that  $\tau$ is a 
 regular
Lyapounov function on $U$.
Then  $\tau$ is still  a  regular
Lyapounov function on $U$  for the closed cone field $\widehat{ \mE_n}$ for $n\geq n_0$.
 There exists a neighborhood $J$ of $t$
and $n_1\geq n_0$ such that $\gamma_n(s)\in U$ for each $s\in J$, $n\geq n_1$.
These properties imply that $\tau\circ \gamma_n$ is non decreasing on $J$
provided $n\geq n_1$. At the limit, we deduce that $\tau\circ \gamma$ is non decreasing on $J$.
\qed

It is useful to control the length of the limit curve:

\begin{lem}\label{lemlength}
	Let $\mC$ be a closed cone field and  $\gamma:[0,1]\lto M$ be a  $\mC$-causal curve which does not contain
	any singular point of $\mC$. There exists $L>0$ such that:

	For each $T\in ]0,1[$, there exists $\epsilon>0$ 
	and an open enlargement $\mE$ of $\mC$ such that 
	each $\mE$-timelike curve $\eta:[0,T]\lto M$ satisfying 
	 $d(\gamma(t), \eta(t))\leq \epsilon $
	for each $t\in [0,T]$ has a length less than $L$.
\end{lem}

\proof
Let $\mE_n$ be a sequence of open cone fields strictly decreasing to $\mC$,
and let $\gamma$ be as in the statement.
We denote  by $\ell(\gamma)$ the length of a curve  $\gamma$.

 We  cover the image of $\gamma$ by finitely many bounded open sets $U_1, \ldots, U_k$ 
 each of which has the property that  there exists a local Lyapounov function $\tau_i$
 on an open neighborhood $V_i$ of $\bar U_i$, which satisfies 
 $
 |v|_x/2\delta_i\geq d(\tau_i)_x \cdot v\geq 2\delta_i |v|_x
 $ 
 for some $\delta_i>0$, 
 and for each $v\in \mC(x), x\in U_i$.
 We set $\delta:=\min \delta_i$ and  prove the statement with $L=(1+\ell(\gamma)/\delta)/\delta$.
 We consider a sequence $\eta_n:[0,1[\lto M$ of $\mE_n$-timelike curves
 converging, uniformly on compact subsets of $[0,1[$, to $\gamma$.
 We have to prove that $\ell(\eta_{n|[0,T]})\leq L$ for $n$ large enough.

 Given  $T\in]0,1[$,
 there exists a finite increasing sequence of times 
 $0=t_0<t_1< \cdots <t_N=T$ such that $\gamma_{[t_j,t_{j+1]}]}$ is contained 
 in one of the open sets $U_i$ for each $j$.
 Then for $n$ large enough, this is also true for $\eta_{n|[t_j,t_{j+1}]}$  and $|\tau_i(\eta(t))-\tau_i(\gamma(t))|\le\frac{1}{2N}$.
We obtain, for $n$ large enough :
 \begin{align*}
 \delta \ell(\eta_{n|[t_j,t_{j+1}]})
 &\leq \tau_i(\eta_n(t_{j+1}))-\tau_i(\eta_n(t_j))
 \leq \tau_i(\gamma(t_{j+1}))-\tau_i(\gamma(t_j))+1/N\\
 &\leq 1/N + \ell(\gamma_{|[t_j,t_{j+1}]})/\delta.
 \end{align*}
 Taking the sum, we obtain that the inequality 
 $$
 \delta \ell(\eta_{n|[0,T]})\leq 1 +\ell(\gamma )/\delta 
 $$
 holds for $n$ large enough, which ends the proof.
 \qed

\begin{prop}\label{proplim}
Let $\gamma_n:[0, a_n[\lto M$ be a sequence of $\mE_n$-timelike curves 
parametrized by arclength, 
such that $\gamma_n(0)$ is bounded and $a_n\lto \infty$.
Then along a subsequence, the sequence $\gamma_n$ converges, uniformly on compact
intervals of $[0,\infty)$, to a limit $\gamma:[0,\infty)\lto M$
which is $\mC$-causal and complete in the sense of Definition \ref{compdef}.
\end{prop}

\proof
Since the curves $\gamma_n$ are $1$-Lipschitz, 
Ascoli Arzela's Theorem  gives, for each $T>0$, the existence of a subsequence
along which $\gamma_n$ converge uniformly on $[0,T]$.
By a diagonal extraction, we get a subsequence along which $\gamma_n$ converge
uniformly on compact intervals.
By Lemma \ref{lem-limit}, the limit $\gamma$ is $\mC$-causal.
Let us prove that this limit is complete.
If it was not complete, it would have finite length and a regular 
limit $\gamma(\infty)=y$ at infinity.
Since the set of regular points is open, there would exit $T>0$ such that 
$\gamma([T,\infty])$ contains only regular points.
We could reparameterize $\gamma$ on $[T,\infty)$ to a curve 
$\tilde \gamma=\gamma\circ \lambda:[0,1[\lto M$,
and extend $\tilde \gamma$  to a causal curve $\tilde \gamma:[0,1]\lto M$.
Lemma \ref{lemlength}, applied  to the causal curve $\tilde \gamma$
and the sequence $\tilde \gamma_n=\gamma_n \circ \lambda$, gives
  $L>0$ such that, for each $S\in ]0,1[$, the curve $\tilde \gamma_{n|[0,S]}$ 
has length
less than $L$ for $n$ large enough.
Observing that $\ell(\tilde \gamma_{n|[0,S]})=\lambda(S)-T$,
this would imply  that $\lambda(S)\leq T+L$ for each $S\in ]0,1[$. 
This is a contradiction since $\lambda$ maps $[0,1[$ onto $[T,\infty)$.
\qed

\begin{cor}\label{cor-complete}
	If $\mC$ is a {non degenerate} closed cone field, then each point is contained in a complete causal curve.
	Each causal curve which is not forward complete can be extended to a forward complete causal curve
	(hence each forward inextendible curve is forward complete).
\end{cor}

\proof
The first point is a consequence of Proposition \ref{proplim} and  Lemma \ref{lem-infinitelengh}.
To prove the second point, we consider a  causal curve $\gamma:[0,T[\lto M$ which is not complete, parametrized by arclength.
Since $\gamma$ has finite length, $T$ is finite.
We consider the limit $y$ of $\gamma$ at $T$, and a forward complete causal curve $\gamma_1$ starting from $y$.
The concatenation of $\gamma$ and $\gamma_1$ is a forward complete causal curve.
\qed

More care is needed in the presence of degenerate points, but we have:

\begin{cor}\label{cor-R}
	Let $\mC$ be a closed cone field.
For each $x\in \mR_{\mC}$, there exists a complete causal curve 
$\gamma$ passing through $x$.
\end{cor}

\proof
Let $\mE_n$ be a sequence of open cone fields 
strictly decreasing to  $\mC$. For each $n$, there exists a 
closed $\mE_n$-timelike curve passing through $x$, that we see as a periodic $\mE_n$-timelike curve 
$\gamma_n:\Rm\lto M$ satisfying $\gamma_n(0)=x$. The curve $\gamma_n$ is periodic and not constant, hence it has infinite length. At the limit, we obtain a complete causal curve passing through $x$.
\qed

The same method also yields:

\begin{cor}\label{cor-back}
Let $Y\subset K$ be two compact sets. If $\mJ_{\mC}^+(Y)$ is contained in the interior of $K$, and 	 $\mF_{\mC}^+(Y)$ is not contained in $K$,
then there exists a backward complete causal curve  $\gamma:]-T,0]\lto M$ contained in $K$ and such that $\gamma(0)\in \partial K$. 
\end{cor}

\proof
For each $n$, there exists an $\mE_n$-timelike curve $\gamma_n:[-T_n,0]\lto K$ such that $\gamma_n(0)\in \partial K$ and $\gamma_n(-T_n)\in Y$, parametrized 
by arclength. If the sequence $T_n$ was bounded, then at the limit we would obtain a $\mC$-causal  curve joining a point of $Y$ to a point of $\partial K$, which 
contradicts the hypothesis that $\mJ_{\mC}^+(Y)$ is contained in the interior of $K$. We deduce that $T_n$ is unbounded, and at the limit we obtain the desired 
backward  complete causal curve.
\qed

 \section{Direct Lyapounov theory}

We consider a closed cone field $\mC$ and explain how to deduce information about stable causality from the existence of appropriate (smooth) Lyapounov 
functions. More precisely we prove the following parts of Theorem~\ref{thm1}, and some variations:
\begin{itemize}
	\item[$\cdot$]
	If there exists a Lyapounov function $\tau$ such that $\tau(x')<\tau(x)$, then $x'\notin \mF^+_{\mC}(x)$.
	\item[$\cdot$]
	If there exists a  Lyapounov function $\tau$ such that $d \tau_x\neq 0$, then $x\notin \mR_{\mC}$.
\end{itemize}

\begin{defn}
An open set $A\subset M$ is a { \it trapping domain} for the open cone field $\mE$ if $\mI^+_{\mE}(A)\subset A$. $A$ is a trapping domain for the closed cone field 
$\mC$ if it is a  trapping domain for some open enlargement $\mE$ of $\mC$.
\end{defn}

In the causality theory of space times such sets are called {\it future sets}, \cite{ms1}.

\begin{lem}\label{L10}
If $A$ is a trapping domain for $\mC$, then there exists an enlargement $\mE$ of $\mC$ such that $ \mI^+_{\mE}(\bar A)\subset A$, in particular,
$\mF_{\mC}^+(\bar A)\subset A$.
\end{lem}

\proof
Let $\mE$ be an open enlargement of $\mC$ such that $\mI^+_{\mE}(A)\subset A$. By Lemma \ref{closure}, $\mI^+_{\mE}(\bar A)= \mI^+_{\mE}(A)\subset A$.
\qed

\begin{lem}\label{lem-trap}
Let $f$ be a $C^1$ function, and $a\in \Rm$. If the inequality $df_x\cdot v>0$ holds for each $x\in f^{-1}(a)$ and $v\in \mC(x)$, then $\{f>a\}$ is a trapping domain.
 	
In particular, if  $a$ is a regular value of the  Lyapounov function $\tau$, then $\{\tau>a\}$ is a trapping domain.
 \end{lem}

The first part of the statement includes the possibility that $\mC(x)$ may be empty, in which case there is no condition on $f$ at $x$.

\proof
Let us consider the open cone field $\mE$ defined
by $\mE(x)=T_xM$ if $f(x)\neq a$ and 
$\mE(x)=\{v\in T_xM: df_x\cdot v>0\}$ if $f(x)=a$.
Our hypothesis on $f$ is that $\mE$ in an  open enlargement of $\mC$.
If $\gamma(t)$ is an $\mE$-timelike curve,
then $f\circ \gamma$ is increasing near each time
$t$ such that $f\circ \gamma(t)=a$. As a consequence, 
if $f\circ \gamma(t)>a$, then $f\circ \gamma(s)>a$
for each $s>t$. This implies that $\{f>a\}$ is a trapping domain.
\qed

\begin{cor}
	Let $\tau$ be a Lyapounov function.
	If $\tau(x')<\tau(x)$, then $x'\not \in \mF^+_{\mC}(x)$.
\end{cor}

\proof
Let $a\in ]\tau(x'), \tau(x)[$ be a regular value of $\tau$ (there exists one by Sard's theorem).
We have $\mF_{\mC}^+(x)\subset \mF_{\mC}^+(\{\tau > a\})\subset \{\tau >a\}$.
\qed

\begin{lem}\label{recreg}
	Let $\tau$ be a Lyapounov function and $x$ a regular point of $\tau$.
	Then there exists a  Lyapounov function $\tilde \tau$ 
	which has the same critical set as $\tau$, 
	and  such that $\tilde \tau (x)$ is a regular value of $\tilde \tau$.
	This implies that $x$ is not stably recurrent.
\end{lem}

\proof
Given a neighborhood $U$ of $x$ on which $\tau$ is regular, let $f$ be a smooth function supported in $U$ and 
such that $f(x)=1$.
For $\delta>0$ small enough, the function $\tau+sf$ is a  Lyapounov function, which is regular on $U$
for each $s\in ]-\delta, \delta[$.
The interval $]\tau(x)-\delta,\tau(x)+ \delta[$ contains a regular value $a$ of $\tau$.
The function $\tilde \tau:= \tau +(a-\tau(x))f$ 
is a Lyapounov function which is regular on $U$.
The number $a:=\tilde \tau(x)$ is a regular value of $\tilde \tau$ :
If $\tilde \tau(y)=a$, then either $y\in U$ and then $d\tilde \tau_y\neq 0$
or $y$ does not belong to the support of $f$, and then $d\tilde \tau_y=d\tau_y\neq 0$
since $a$ is a regular  value of $\tau$.
\qed

\section{Smoothing} 

The goal of the present section is to prove the following regularization statement, which is one of our  main 
technical tools to prove the existence of Lyapounov functions. We work with a 
closed cone field $\mC$ on 
the manifold $M$. We say that the open set $A\subset M$ is smooth 
if its boundary is a smooth submanifold, we say that the open set $A\subset M$ is smooth 
near the set $X$ if there exists 
an open set
$U$ containing $X$ such that $U\cap A$ is smooth in $U$.

\begin{prop}\label{prop-Areg}
	Let $A_0$ be a trapping domain, let $F_i$ be a closed set contained in $A_0$, let $F_e$
	be a closed set disjoint from $\bar A_0$, and let $\theta_0$ be a point in the boundary of $A_0$.
	
	Then there exists a smooth (near $\mD(\mC)$) trapping domain $A'_0$ which contains $F_i$, whose boundary contains 
	$\theta_0$, and whose closure is disjoint from $F_e$.
\end{prop}

The proof is given in  \ref{sec-proof-prop-Areg} after the exposition of preliminary material.
 
\subsection{Local properties of trapping domains}

Given an open set $A\subset M$ and $x\in \partial A$, we say that $A$ is {\it locally trapping at $x$} if 
one of the  following  equivalent conditions hold:

\begin{itemize}
\item[(i)] There exists an open cone field $\mE$ which contains $\{x\}\times \mC(x)$ and such that $\mI_{\mE}^+(A)\subset A$.
\item[(ii)] There exists a compact neighborhood $K$ of $x$ such that $A$ is trapping for $\mC_K$ (the cone field equal to $\mC$
on $K$ and degenerate outside of $K$).
\item[(iii)] There exists an neighborhood $U$ of $x$ such that $A\cap U$ is trapping for the cone field $\mC$ on $U$.
\end{itemize}

\begin{lem}\label{lem-loc}
The open set $A$ is a trapping domain   for $\mC$  if and only if 
it is locally trapping at each point $x\in \partial A$.
\end{lem}

 \proof
 If $A$ is a trapping domain, then there exists an open enlargement $\mE$ of $\mC$ such that $\mI^+_{\mE}(  A)\subset A$.
This implies that $\mC$ is locally trapping $A$ at each point of $\partial A$.
 
 Let us now prove the converse. For each point $x\in \partial A$, there exists an  open cone field $\mE_x$ such that $\mC(x)\subset \mE_x(x)$ and $\mI^+_{ \mE_x}(A)\subset A$.
The inclusion $\mC(y)\subset \mE_x(y)$ then holds for all $y$ in an open neighborhood $U_x$ of $x$ in $\partial A$.
 We consider a sequence $x_i$ such that  the open sets $U_{x_i}$ form a locally finite covering of $\partial A$.
 For each $x\in \partial A$, we denote by $J(x)$ the finite set of indices such that $x\in \bar U_{x_i}$.
 Since the covering is locally finite, there exists a neighborhood $V$ of $x$ in $\partial A$ which is disjoint from $U_{x_i}$
 for each $i\notin J(x)$.
 We define, for each $x\in \partial A$, the open cone $\mE(x):= \bigcap _{i\in J(x)} \mE_{x_i}(x)$. For $x\notin \partial A$, we set $\mE(x)=T_xM$.
 We claim that $\mE:= \bigcup_{x\in M} \{x\}\times \mE(x)$ is an open cone field.
 Indeed, for each $x\in A$, the intersection 
 $ \bigcap_{i\in J(x)} \mE_{x_i}$ is an open cone field which is contained in $\mE$
 in a neighborhood of $x$, and equal to $\mE$ at $x$.
 
 By construction, $\mE$ is an enlargement of $\mC$. Let us verify that $\mI^+_{\mE}(A)\subset A$.
 If not, there exists an $\mE$-timelike curve $\gamma$ such that $\gamma(t)\in A$ on $[0,T[$ 
 and $\gamma(T) \in \partial A$.
 We have $\dot \gamma(T)\in \mE(\gamma(T))\subset \mE_{x_i}(\gamma(T)))$
 for some  $i$ (any $i$ such that $\gamma(T)\in U_{x_i})$.
 For this fixed $i$, the curve $t\mapsto\gamma(t)$ is then $\mE_{x_i}$-timelike
 on $[S,T[$ for some $S<T$.
 This contradicts the inclusion $\mI^+_{\mE_{x_i}}(A)\subset A$.
 \qed
 
 Let  $\mE$ be an non degenerate open cone field,
 and $A$ be a trapping  domain for $\hat \mE$.  Then $\hat\mE(x)$ is regular for all $x\in \partial A$. By Lemma
 \ref{localcone}, at each point $x\in \partial A$, 
 there exists a chart 
  $\phi: B^{d-1}(2)\times ]-2,2[\lto M$ which sends $(0,0)$ to $x$ and has the property that
 $$
 Q_1\subset   \phi^* \hat \mE(y,z)\subset  Q_0
 $$
 for all $(y,z)\in B^{d-1}(2)\times ]-2,2[$. We recall that $B^d(r)$ is the open ball 
 of radius $r$ centered at $0$ in $\Rm^d$ and that 
  $Q_s, s\geq 0$  is the  open cone $Q_s=\{(y,z)\in \Rm^{d-1}\times \Rm :z> s|y|\}\subset \Rm^d$.

 \begin{lem}\label{lembox}
 	There exists a $1$-Lipschitz function $g: B^{d-1}(2)\lto ]-2,2[$ such that 
 	$\phi^{-1}(A)$ is the open epigraph $\{z>g(y)\}$, hence $\phi^{-1}(\partial A)$ is the graph of $g$.
 	Note that $g(0)=0$.
 \end{lem}

 \proof
 Let us define the function  $g(y)= \inf \{ z\in ]-2,2[: \phi(y,z)\in A\}$. Since $\phi(0,0)\in \partial A$ and $Q_1 \subset  \phi^* \hat{\mE}$ there exists for all $y\in B^{d-1}(2)$ a $z<2$ with 
 $\phi(y,z)\in A$.
 Then $\phi(y,g(y))\in \bar A$ for each $y\in B^{d-1}(2)$.
 The curve $t\lmto \phi(y,z+t)$ is $\bar \mE$-causal 
 hence the set $\{(y,z): g(y)<z<2 \}$ is contained in $A$.
 Furthermore, since $Q_1\subset  \phi^* \bar \mE$, 
 the curve $\phi(y+tv,g(y)+t)$ is  $\bar \mE$-causal for each $y \in B^{d-1}(1)$ and $v\in \bar B^{d-1}(1)$.
 This implies that $g$ is $1$-Lipschitz.
 \qed

Let $\mC$ be a closed cone field on $\Rm^{d-1}\times \Rm$, and let 
$A\subset \Rm^{d-1}\times \Rm$ be a trapping domain which is 
the open epigraph of the Lipschitz function
$g:\Rm^{d-1}\lto \Rm$.

\begin{lem}\label{lem-inC}
	For each point $x=(y,g(y))$ of $\partial A$, the following statements are equivalent:
	\begin{itemize}
		\item[(a)]
		The domain $A$ is locally trapping at $x$.
		\item[(b)]
		The inequality $v_z >p\cdot v_y$ holds for each  $(v_y,v_z)\in \mC(x)$,
		and each $p\in \partial g(y)$.
	\end{itemize}
\end{lem}

\proof
If (b) does not hold, there exists $p\in \partial g(y)$ and $(v_y,v_z)\in \mC(x)$ such that $p\cdot v_y \geq v_z$.
Then for  each open cone field $\mE$ containing $\{x\}\times \mC(x)$, there exists $w=(w_y,w_z)\in \mE(x)$ such that $w_z<p\cdot w_y$.
Then there exists an open interval $I$ containing $0$ and a neighborhood $V$ of $x$ such that the curve $t\lmto x'+tw$ is $\mE$-timelike 
on $I$ for each $x'\in V$. The inequality $w_z<p\cdot w_y$ implies that the Clarke differential of the function $t\lmto g(y+tw_y)-tw_z$ contains 
a positive value. As a consequence, this function is not non increasing in any neighborhood of $0$. In other words,  there exists $t_1<t_2$ in $I$ such
that 
$$g(y+t_1w_y)-t_1w_z<g(y+t_2w_y)-t_2w_z.$$
We can assume moreover that $t_1$ is sufficiently small to have 
 $(y,g(y+t_1w_y)-t_1w_z)\in V$.
This implies that the curve 
$$\eta=(\eta_y,\eta_z): t\lmto (y+tw_y,g(y+t_1w_y)+(t-t_1)w_z)$$
is $\mE$-timelike on $I$. We observe  that $\eta_z(t_1)=g(\eta_y(t_1))$ and $\eta_z(t_2)<g(\eta_y(t_2))$.
As a consequence, we do not have $\mI^+_{\mE}(\bar A)\subset A$. We have proved that $A$ is not locally trapping at $x$.

Conversely, let us assume that (b) holds and consider the  cone 
$$\Omega=\{(v_y,v_z): v_z> p\cdot v_y, \forall p\in \partial g(y)\}.$$
Since $\partial g(y)$ is compact, this is an open cone.
We consider an open cone $\Omega_1$ such that 
$$\mC(x)\subset  \Omega_1\subset  \hat\Omega_1\subset  \Omega.$$
In view of the semi-continuity of the Clarke differential,
there exists an open  neighborhood $U$ of $y$ in $\Rm^{d-1}$ such that 
the inequality $v_z>\sup _{p\in \partial g(y')}p\cdot v_y$ 
holds for each $(v_y,v_z)\in \Omega_1$ and each $ y'\in U$.
We consider the open cone field $\mE$ which is equal to $\Omega_1$ on $U\times \Rm$
and empty outside, and prove that $\mI^+_{\mE}(A)\subset A$.
Otherwise, there exists a smooth curve  $\gamma=(\gamma_y,\gamma_z)$, which is
timelike  for $\mE$, 
and such that $\gamma_z(T)= g(\gamma_y(T))$ and $\gamma_z(t)> g(\gamma_y(t))$
for each $t\in [0,T[$.
Then, we have  $\gamma_y(T)\in U$ and $(\dot \gamma_y(T), \dot \gamma_z(T))\in \Omega_1$
hence 
$$\dot \gamma_z(T) > \sup_{p\in \partial g(\gamma_y(T))}p\cdot \dot \gamma_y(T).$$
This implies that the function $\gamma_z(t)-g(\gamma_y(t))$ is increasing
near $t=T$,  a contradiction.
\qed

\subsection{De Rham Smoothing}\label{secreg}

\begin{prop}\label{prop-rham}
	For each Lipschitz function $g:\Rm^ d\lto \Rm$, there exists a family $g_s, s>0$ of Lipschitz functions on $\Rm^d $
	which converge uniformly to $g$ as $s\lto 0$ and such that:
	\begin{itemize}
		\item[(a)]
		$g_s$ is smooth on $B^d(1)$ and equal to $g$ outside of this ball  for each $s>0$, and moreover $g_s$ is smooth on any open subset $O\subset \Rm^d$ where $g$ 
		is already smooth.
		\item[(b)]
		$\limsup _{s\lto 0}(\textnormal{Lip } g_s )\leq \textnormal{Lip }g$.
		\item[(c)]
		If $V\subset \Rm^d\times (\Rm^d)^*$ is an open set containing the graph
		 $\partial g:=\{(x,p): x\in \bar B^{d}(1), p\in \partial g(x)\}$ of the Clarke differential of $g$,
		then $V$ contains the graph
		 $\partial g_s:=\{(x,p): x\in \bar B^{d}(1), p\in \partial g_s(x)\}$ for $s$ small enough.
	\end{itemize}
	If $y_1,\ldots,y_N$ are finitely many points in $B^d(1)$, then we can assume in addition
	that $g_s(y_i)=g(y_i)$ for each $i=1,\ldots, N$ and each $s>0$.
\end{prop}

\proof
We use { the de} Rham smoothing procedure.
We follow the notations of { \cite[Lemma A.1]{BZ}}.
There exists a smooth  action { $a\colon \Rm^d\times \Rm^d\to \Rm^d$, $(y,x)\mapsto a(y,x)$} of $\Rm^d$ on itself 
(meaning that $a(y,a(y',x))=a(y+y',x)$)
 such that:
\begin{itemize}
	\item[$\cdot$]
	$a(y,x)=x$ for each $y\in \Rm^d$ and $x \in \Rm^d-B^d(1)$
	\item[$\cdot$]
	The action of $\Rm^d$ on $B^d(1)$ is  conjugated to the standard 
	action of $\Rm^d$ on itself by translations
	(there exists a diffeomorphism $\varphi:B^d(1)\lto \Rm^d$ such that
	$\varphi \circ a(y,\varphi(x))=y+\varphi(x)$).
	\item[$\cdot$]  The diffeomorphisms {$a_y\colon x\to a(y,x)$} converge to the identity  $C^1$-uniformly  for $y\lto 0$.
\end{itemize}

Given a Lipschitz function $g:\Rm^d\lto \Rm$,
we define 
$$
g_s(x):= \int_{\Rm^d} s^{-d} g(a(y,x))\rho(-y/s) dy
$$
where $\rho$ is a mollification kernel supported in $B^d(1)$.
Properties (a) and (b) are proved, for example, 
in  {\cite[Lemma A.1]{BZ}}.
Let us now prove property (c).

We denote by $V(x)\subset (\Rm^d)^*$ the set of $p$ such that $(x,p)\in V$.
We cover the compact set $\bar B^d(1)$ by finitely many balls $B_i$ each of which has the following property:
There exists  convex  open sets $W_i$ and $V_i$ in $ (\Rm^d)^*$ such that  $\partial g(x)\subset W_i\subset \bar W_i\subset V_i\subset V(x)$
for each $x\in 2B_i$ (the ball of same center and double radius).

For each $i$, we define $n_i(v):= \sup_{p\in W_i} p\cdot v$ and $m_i(v)= \sup_{p\in V_i} p\cdot v$
which are convex and positively one-homogeneous (hence subadditive) functions. { Since
$\overline{W_i}\subset V_i$ there} exists $\delta>0$ such that $m_i(v)\geq n_i(v)+\delta|v|$.
Note that $V_i$ (resp. $W_i$) is precisely the set of linear forms $p$ satisfying $p\cdot v \leq m_i(v)$ (resp. $n_i(v)$) for each $v$ .
The function $g$ is $n_i$-Lipschitz on $2B_i$, which means that 
$$
g(x')-g(x)\leq n_i(x'-x)
$$
for each $x$ and $x'$ in $2B_i$.
Since the diffeomorphisms $a_y$ converge to the identity  $C^1$-uniformly  as $y\lto 0$,
we have
$$
\big|a(y,x')-a(y,x)-x'+x\big|=\left|\int _0^1 (\partial_x a(y,x+t(x'-x))-Id) \cdot (x'-x) dt\right|
\leq {\epsilon}(|y|)|x'-x|
$$
with a function ${\epsilon}$ converging to $0$ at $0$.
For $s$ small enough, we have $a(y,x)\in 2B_i$ for each $x\in B_i$ and $|y|\leq s$, and ${\epsilon}(s)<\delta$.
We then obtain, for $x$ and $x'$ in $B_i$,
\begin{align*}
|g_s(x')-g_s(x)|
&\leq \int s^{-d}\big| g\circ a_y (x')- g\circ a_y (x)\big| \rho(-y/s) dy\\
&\leq \int s^{-d} n_i\big(a_y (x')- a_y (x)\big) \rho(-y/s) dy\\
&\leq \int s^{-d} m_i(x'-x) \rho(-y/s) dy=m_i(x'-x).
\end{align*}
This implies that $dg_s(x)\cdot v \leq m_i(v)$ for each $v$, at each point of differentiability $x$ of $g_s$ in $B_i$, hence that 
$\partial g_s(x)\subset V_i\subset V(x)$ for each $x\in B_i$.
Since the covering $B_i$ is finite, this inclusion holds for all $x\in \Rm^d$ provided $s$ is small enough.

The function $g_s$ constructed so far does not necessarily satisfy the additional conditions $g_s(y_i)=g(y_i)$.
We thus consider the modified function
$$\tilde g_s(x)=g_s(x)+\sum_{i=1}^N(g(y_i)-g_s(y_i))h_i(x),
$$
where $h_i, 1\leq i\leq N$ are  non negative smooth {functions} supported on $B^d(1)$ and
satisfying $h_i(y_i)=1$ and $h_i(y_j)=0$ for $j\neq i$.
This modified family of functions satisfies the three points of the statements
since $g_s(y_i)\lto g(y_i)$ for each $i$, and
$$
\partial \tilde g_s(x)=\partial g_s(x)+\sum_i (g(y_i)-g_s(y_i))dh_x
$$
for each $x$.
\qed

\subsection{Proof of Proposition \ref{prop-Areg}}\label{sec-proof-prop-Areg}

We first give the proof under the assumption that $\mD(\mC)=M$.
Since $A_0$ is also a trapping domain for some open enlargement $\mE$
of $\mC$  (Lemma \ref{L10}), we can assume without loss of generality that $\mC$
is the closed hull of a non degenerate  open cone field.
 
We consider a locally finite covering of $\partial A_0$ by  domains 
$$U_k(1)=\phi_k(B^{d-1}(1)\times ]-1,1[)$$
associated to charts $\phi_k: B^{d-1}(2)\times ]-2,2[\lto M, k\geq 1$ which have the property that
$$
Q_1\subset   \phi_k^*  \mC(y,z)\subset  Q_0
$$
for all $(y,z)\in B^{d-1}(2)\times ]-2,2[$.
We denote by $x_k$ the points $\phi_k(0)$, $k\geq 1$ and set $x_0=\theta_0$.
We moreover assume that the open sets $U_k(2):=\phi_k(B^{d-1}(2)\times ]-2,2[)$
are all disjoint from $F_i$ and $F_e$.

By Lemma \ref{lembox},
the open set $\phi_1^{-1}(A_0)$ is the epigraph of a $1$-Lipschitz function 
$f_1:B^{d-1}(2)\lto ]-2,2[$ such that $f_1(0)=0$.
The bounded set $U_1(1)$ contains finitely many of the points $x_i$.
We denote by $y_1, \ldots,y_N$ the first component of the preimages of these points.
So those of the points $x_i$ which are contained in $U_1(1)$ are 
$\phi_1(y_1,f_1(y_1)),\ldots ,\phi_1(y_N,f_1(y_N))$.

For each $y\in B^{d-1}(2)$, the inequality
$v_z>p\cdot v_y$ holds for each $(v_y,v_z) \in \phi_1^*\mC(y,f_1(y))$ and each 
$p\in \partial f_1(y)$ by Lemma \ref{lem-inC}.
By a compactness argument, we find an open neighborhood $W_1$ of
 $\{(y,f_1(y)), y\in \bar B^{d-1}(1)\}$,
 and an open neighborhood $V_1$ of 
 $$\partial f_1:= \{(y,p), y\in \bar B^{d-1}(1), p\in \partial f_1(y) \}$$
 with the property that the inequality
$
v_z> p\cdot v_y
$
holds for each $x=(y,z)\in W_1$, each $p\in V_1(y)$, and each 
$(v_y,v_z)\in \phi_1^*\mC(x)$. We have denoted by $V_1(y)$ the set of linear forms 
$p$ such that $(y,p)\in V_1$.

By Proposition \ref{prop-rham}, there exists a function $g_1:B^{d-1}(2)\lto \Rm$
which is $2$-Lipschitz, smooth on $B^{d-1}(1)$, equal to $f_1$ outside of 
$B^{d-1}(1)$, and 
satisfies:
\begin{itemize}
	\item[$\cdot$] $\partial g_1=\{(y,p), y\in \bar B^{d-1}(1) , p\in \partial g_1(y)\} \subset V_1$, 
	\item[$\cdot$] $(y,g_1(y))\in W_1$ for each $y\in B^{d-1}(1)$,
	\item[$\cdot$] $g_1(y_j)=f_1(y_j)$ for $j=0, \ldots, N$.
\end{itemize}

In particular, $g_1(0)=0$, hence $g_1$ takes vales in $]-2,2[$.

Let $A_1$ be the open set such that $A_1\cap (M-U_1(1))=A_0\cap (M-U_1(1))$
and such that $\phi_1^{-1}(A_1)$ is the open epigraph of $g_1$.
The domain $A_1$ is locally trapping at each point of its boundary.
Indeed, such a point $x$ either belongs to $\partial A_0\cap (M-\bar U_1(1))$, and then $A_1=A_0$
near $x$, or it is of the form $\phi_1(y,g_1(y))$ for some  $y\in \bar B^{d-1}(1)$.
In this second case, we have $\phi_1^{-1}(x)=(y,g_1(y))\in W_1$,
hence the inequality $v_z>p \cdot v_y$ holds for each $(v_y,v_z) \in \phi_1^*\mC(x)$
and each $p\in \partial g_1(y) \subset V_1(y)$.
The conclusion then follows from  Lemma \ref{lem-inC}.
We deduce by Lemma \ref{lem-loc} that $A_1$ is a  trapping domain for $\mC$.

By the same method, we build inductively a sequence $A_m, m\geq 0$ of trapping domains 
which have the following properties:
\begin{itemize}
\item[$\cdot$] $\partial A_m$ contains all the points $x_k$ (hence the point $\theta_0$), $F_i$ is contained in $A_m$ and $F_e$ is disjoint from $\bar A_m$.
\item[$\cdot$] The boundary $\partial A_m$ is contained 
in $\bigcup_{k\geq 1} U_k(1)$, and its intersection with $\bigcup_{m\geq k\geq 1} U_k(1)$
 is a smooth hypersurface.
\item[$\cdot$] The symmetric difference between $A_m$ and $A_{m-1}$ is contained in $U_m(1)$.
\end{itemize}

We denote by $A'_{0}:=\liminf A_m$ the set of points $x$ which belong to all but finitely many
of the sets $A_m$.
We claim that $A'_{0}$ satisfies the conclusions of Proposition \ref{prop-Areg}.
Since the covering $U_k(1)$ is locally finite, the intersection $A_m\cap K$ stabilizes to $A'_{0}\cap K$
for each compact $K$, \textit{i. e.} $K\cap A_m= K\cap A'_{0}$ for all $m$ large enough.
This implies that $A'_{0}$ is open, and that $\partial (A'_0)=\liminf \partial (A_m)$.
This boundary is smooth, contains all the points $x_k$, and is contained in $\bigcup_{k\geq 1} U_k(1)$.

To prove that $A'_{0}$ is a trapping  domain, it is enough to observe that it is locally
trapping at each point $x$ of its boundary.
Since the sequence $A_k$ stabilizes in a neighborhood of $x$,
this follows from the fact that each of the open sets $A_k$ is trapping.
This ends the proof of Proposition \ref{prop-Areg} under the assumption that $\mD(\mC)=M$.

In case that $\mD(\mC)\neq M$ we consider an enlargement $\mE$ of $\mC$
such that $A_0$ is a trapping domain for $\hat \mE$  (Lemma \ref{L10} and Lemma \ref{intermediate}).
We can apply the result just proved  on the manifold $\mD(\mE)$, to the cone field $\hat \mE$.
We deduce the existence of a smooth trapping region $A'_0$ for $\hat \mE$ in $\mD(\mE)$
which contains $F_i\cap \mD(\mE)$, is disjoint from $F_e \cap \mD(\mE)$, and whose boundary contains $\theta_0$.
Let $O$ be an open subset of $M$ which contains $F_i $ and whose closure is disjoint from $F_e$,
and let $Z\subset \mD(\mE)$ be a closed neighborhood of $\mD(\mC)$ in $M$.
The open set $A'_0\cup ((M-Z)\cap O)$ then satisfies the conclusions of   
Proposition \ref{prop-Areg}.
\qed

\section{Existence of Lyapounov functions}\label{sec-ex}

We consider in this section a closed cone field $\mC$ and prove several existence results for Lyapounov functions,  in particular Theorems 
\ref{thm1}, \ref{thm2} and \ref{thm3}.

\subsection{Smooth trapping domains and Lyapounov functions}

We associate (smooth)  Lyapounov functions to   smooth trapping domains:

\begin{prop}\label{prop-fol}
	Let $A$ be  smooth trapping domain, then 
	there exists a (smooth) Lyapounov function $\tau:M\lto [-1,1]$
	such that $A=\{\tau >0\}$ and
	all values in $]-1,1[$ are 
	regular values of $A$
	$($hence $\partial A=\{\tau=0\})$.
	
	If $F_i$ and $F_e$ are closed sets  
	contained in $A$ and disjoint from $\bar A$, respectively, we can moreover impose that $\tau=1$ on $F_i$ and 
	$\tau=-1$ on $F_e$.
\end{prop}

\proof
We consider a collar of the hypersurface $H:=\partial A$ in the manifold $M-(F_e\cup F_i)$, that is a smooth embedding
$\psi: H\times \Rm\lto M-(F_e\cup F_i)$ such that $\psi(H\times \{0\})=\partial A$
and $\psi^{-1}(A)=H\times ]0, \infty)$.
We will prove the existence of a Lyapounov function 
$\tilde \tau:H\times \Rm\lto [-1,1]$ for the cone field $\psi^* \mC$, which has the following properties:
\begin{itemize}
	\item[$\cdot$] $\tilde \tau =0$ on $H\times \{0\}$.
	\item[$\cdot$]
	$\tilde \tau=1$ on $H\times [1,\infty)$ and $\tilde \tau = -1$ on $H\times (-\infty, -1]$.
	\item[$\cdot$] The values in $]-1,1[$ are regular for $\tilde \tau$.
\end{itemize}
Assuming the existence of the function $\tilde \tau$, we obtain the Lyapounov function $\tau$
on $M$ as follows:
$ \tau=\tilde \tau \circ \psi^{-1}$ on $U=\psi(H\times \Rm)$, $\tau=1$ on
$A-U$, and $\tau=-1$ on $M-(A\cup U)$. 

Let us now prove the existence of the Lyapounov function $\tilde \tau$ on $H\times \Rm$.
We denote by $(y,z)$ the points of $H\times \Rm$.
The cone field  
$$\tilde \mC(y,z)=\psi^*\mC(y,z)=(d\psi_{(y,z)}^{-1}\cdot \mC(\psi(y,z)))
$$
is a closed cone field on $H\times \Rm$.
The cones $\tilde \mC(y,0)$ satisfy $v_z>0$ for each $(v_y,v_z)\subset \tilde \mC(y,0)$.
Fixing a Riemannian metric on $H$, 
there exists a smooth positive function $\delta(y)$ on {$H$} such that 
$$\tilde \mC(y,0)\subset \{(v_y,v_z): v_z\geq 3\delta(y) \|v_y\|\}
$$
for each $y\in H$.
Then, there exists a smooth positive function {$\epsilon$ on $H$} such that
$$\tilde \mC(y,z)\subset \{(v_y,v_z): v_z\geq 2\delta(y) \|v_y\|\}
$$
provided $|z|\leq \epsilon(y)$.
Let $f:H\lto \Rm$ be a smooth positive function such that 
$\|df_y\|\leq \delta(y)$ and $f(y)\leq \epsilon(y)$ for all {$y\in H$},
see Lemma \ref{lemsmall} below for the existence of such a function.
We set 
$$\tilde \tau(y,z)=\phi(z/f(y)),
$$
where $\phi:\Rm\lto[-1,1]$ is a smooth nondecreasing function which has positive derivative on $]-1,1[$
and is equal to $1$ on $[1,\infty)$ and to $-1$ on $(-\infty,-1]$.
The set of regular points of the function $\tilde \tau$ is $\{(y,z): |z|<f(y)\}=\tilde \tau ^{-1}(]-1,1[)$.
At such a point $(y,z)$, we compute
$$d\tilde \tau_ {(y,z)}\cdot(v_x,v_z)=\frac{\phi'(z/f(y))}{f(y)}
\big(
v_z-\frac{z}{f(y)}df_y\cdot v_y
\big)\geq \frac{\phi'(z/f(y))}{2f(y)}v_z
$$ 
for $(v_y,v_z)\in \tilde \mC(y,z)$ since
$
|(z/f(y))df_y\cdot v_y|\leq \delta(y) \|v_y\|\leq v_z/2.
$
\qed

We will also need a variant of the above result.

\begin{cor}\label{lyapcor}
	Let $A$ be a trapping domain which is smooth near $\mD(\mC)$.
	Then there exists a Lyapounov function $\tau:M\lto [-1,1]$ such that $A=\{\tau>0\}$ and such that
	$\tau$ is regular at each point of $\tau^{-1}(]-1,1[)\cap \mD(\mC)$.
	
	If $F_i$ and $F_e$ are closed sets 
	contained in $A$ and disjoint from $\bar A$, respectively, we can moreover impose that $\tau=1$ on $F_i$ and 
	$\tau=-1$ on $F_e$.
\end{cor}

We recall the classical:

\begin{lem}\label{L11}
For $A\subset M$ open there exists a non negative smooth function $f\colon M\to \Rm$ 
such that $A=\{f>0\}$.
\end{lem}

\proof
Choose a locally finite open cover $\{B_i\}_i$ of $A$ and
a subordinate partition of unity $\{\lambda_i\}_i$. There exists a positive sequence $\{a_i\}_i$
(see \cite{fathi97} for example) such that 
$f:=\sum_i a_i \lambda_i$
is smooth. It satisfies the desired property.
\qed

\textsc{Proof of Corollary \ref{lyapcor}.}
Let $U$ be an open neighborhood of $\mD(\mC)$ such that $\partial A\cap U$ is smooth.
Let $V$ be the complement of $\mD(\mC)$.
Let $T:U\lto [-1,1]$ be a  Lyapounov function of $\mC$ on $U$ such that 
$A\cap U=\{T>0\}$ and all points in $]-1,1[$ are regular values of $T$.
 We obtain such a function by applying  Proposition \ref{prop-fol} on the manifold $U$.
{Using Lemma \ref{L11} we can choose a smooth function $f$} on $M$ such that $f=1$ on $F_i$, $f>0$ on $A$, $f<0$ outside of $\bar A$,
and $f=-1$ on $F_e$.
Let $g,h$ be a partition of  unity associated to the open covering $(U,V)$ of $M$.
We set $\tau=g T+h f$.
\qed

\begin{lem}
	\label{lemsmall}
	
	Le $N$ be a Riemannian manifold, and let $\epsilon(x)$ be a positive continuous function
	on $N$. There exists a smooth function $f:N\lto \Rm$ such that $|f(x)|\leq \epsilon(x)$
	and $|df(x)|\leq \epsilon(x)$ for each $x$.
\end{lem}

\proof
We consider a locally finite partition of the unity by smooth 
compactly supported function $g_i$, and  
a function $f$ of the form 
$f=\sum_i a_i g_i$ for some positive sequence $a_i$.  

We set $h(x):= 1+\sum_i |dg_i(x)|$. This is a continuous positive function 
(the sum is locally finite). We claim that the function $f=\sum a_ig_i$  
satisfies the desired inequalities provided $0<a_i < \min_{x\in K_i} \epsilon/h$,
where $K_i$ is the support of $g_i$. 

To prove the claim, we define for each $x$ the finite set $I(x)$ of indices $i$ such that $x\in K_i$.
For each $i\in I(x)$, we have $a_i < \epsilon(x)/h(x)$, hence
$$f(x)=\sum_{i\in I(x)}a_ig_i(x) < (\epsilon(x)/h(x))\sum _i g_i(x)< \epsilon(x).$$
Moreover, $|df(x)|\leq  \sum_{i\in I(x)}a_i|dg_i(x)|<\epsilon(x)$.
\qed

\subsection{Conley Theory for closed cone fields}\label{subsec-ex}

We prove Theorem \ref{thm1} and \ref{thm2}.
We say that $a$ is a {{\it relative regular value}} of $\tau$ if $\tau^{-1} (a)\cap \mD(\mC)$
consists of regular points of $\tau$.

\begin{prop}
If $x$ is not stably recurrent, then there exists a Lyapounov function $\tau$ such that $\tau(x)$ is
in the interior of the set of relative regular values of $\tau$ (in particular, $\tau$ is regular at $x$).
\end{prop}

\proof
There are two cases.
 Either $\mC(x)$ is degenerate, or there exists an enlargement $\mE$ of $\mC$ such that $x\not\in \mI^+_{\mE}(x)$
 and such that $\mE(x) \neq \emptyset$.
 
 In the first case, the point $x$ belongs to the open set  $M-\mD(\mC)$.
 Then there exists a smooth function $\tau$ compactly supported inside this open set, 
 and such that $\tau(x)$ is in the interior 
 of the set of regular values of $\tau$. This function $\tau$ is a Lyapounov function.
 
 In the second case, the  set  $A_0:=\mI^+_{\mE}(x)$
is a trapping domain for $\mC$ whose boundary contains $x$.
Proposition \ref{prop-Areg} gives the existence of a 
trapping domain which is smooth near  $\mD(\mC)$ and  whose boundary contains $x$.
Corollary \ref{lyapcor} then
implies the existence of the desired Lyapounov function.
\qed

\begin{prop}\label{propstep}
Let $x$ and $x'$ be two points such that $x'$ does not belong to $\mF_{\mC}^+(x)$. Then there exists a Lyapounov function $\tau:M\lto [-1,1]$
such that $\tau(x')=-1$, $\tau(x)=1$, and all values in $]-1,1[$ are relative regular values of $\tau$.
\end{prop}
\proof
We consider  two cases. Either $\mC(x)=\emptyset$ or there exists
an enlargement $\mE$ of $\mC$ such that  $x'\notin \mI_{\mE}^+(x)\cup \{x\}$ and $\mE(x)\neq \emptyset$.
 
 In the first case, we take a smooth function $\tau$ 
 which is equal to $1$ in a small neighborhood of $x$
 and  $-1$ in a neighborhood of $\mD(\mC)\cup \{x'\}$.

 In the second case, the set  $A_0:=\mI^+_{\mE}(x)$ is a trapping domain containing $x$ in its closure and not containing $x'$.
Proposition \ref{prop-Areg} then  implies the existence of a smooth 
(near $\mD(C)$) trapping domain containing $x$ in its closure and not containing $x'$.
Corollary  \ref{lyapcor}
implies the existence of a   Lyapounov function $\tilde \tau:M\lto [-1,1]$ such that
$\tilde \tau(x')\leq 0$ and   $\tilde \tau(x)\geq 0$ and values in $]-1,1[$ are relatively regular.
By slightly perturbing $\tilde \tau$ near $x$ if necessary, we can assume that 
$\tilde \tau(x)>0$. We then set $\tau= f\circ \tilde \tau$, with a non-decreasing 
smooth function $f: \Rm \lto [-1,1]$ which has positive derivative on
$]\tilde \tau(x'), \tilde \tau(x)[$ and sends this interval onto 
$]-1,1[$.
\qed

Theorem \ref{thm1} obviously follows from the two propositions above.
Let us prove Theorem \ref{thm2}.

\textsc{Proof of Theorem \ref{thm2}.}
Let us consider the set $\mL$ of  Lyapounov functions $\tau$ which have the property that 
they take values in $[-1,1]$ and that each point of $\tau^{-1}(]-1,1[)\cap \mD(\mC)$
is regular for $\tau$ (in other words, values in $]-1,1[$ are relative regular values).
We endow $\mL$ with the topology of $C^1$ convergence on compact sets. 
Being a subset of the separable metric space $C^1_{loc}(M,\Rm)$, it is a separable metric space.
We consider a dense sequence $\tau_i$ in $\mL$.
There exists a positive sequence $a_i$ such that $\tau=\sum a_i \tau_i$ converges in $C^k$ for each $k$ on each compact set
(see \cite{fathi97} for example).
We can moreover assume that $a_{i+1}\leq a_i/5$. We claim that the sum $\tau$ then satisfies all 
the conclusions of Theorem \ref{thm2}.

For each point $x\in \mD(\mC)$ which is not stably recurrent, there exists a Lyapounov function $f\in \mL$ such that $df_x\neq 0$, by Corollary \ref{lyapcor}.
As a consequence, there exists $i$ such that $d\tau_i(x)\neq 0$.
If $v\in \mC(x)$ is  not zero, then all terms of the sum  $d\tau_x\cdot v=\sum_i a_id\tau_i(x)\cdot v$ 
are non negative, and one of them is positive, hence the sum is positive. We deduce that $x$ is a regular point of $\tau$.

Let us then consider two points $x\neq x'$ in $M$ such
that $x'\in \mF^+(x)$ and $x\notin \mF^+(x')$. The first point implies that
$\tau(x')\geq \tau(x)$ for each Lyapounov function $\tau$.
The second point implies, by Proposition \ref{propstep},
the existence of a function $f\in \mL$ such that $f(x')=1, f(x)=-1$. 
By density,  there exists  $j$ such that $\tau_j(x')>\tau_j(x)$.
The difference $\tau(x')-\tau(x)$ is thus the sum of non negative terms
one of which is positive, hence $\tau(x')>\tau(x)$.

Finally, if $x$ and $x'$ are two stably recurrent points which do not belong to the same stable class,
then there exists $i$ such that $\tau_i(x)\neq \tau_i(x')$.
We consider the first index $j$ with this property.
Since $x'$ and $x$ are stably recurrent, we necessarily have 
that $\tau_j(x')=\pm 1$ and $\tau_j(x)=\pm 1$ (the only critical values of $\tau_j)$.
We assume for definiteness that $\tau_j(x')=1, \tau_j(x)=-1$.
Then
$$
\tau(x')-\tau(x)=\sum_i a_i (\tau_i(x')-\tau_i(x))\geq a_j-\sum _{i>j} a_i\geq 3a_j/4>0
$$
since $a_i \leq a_j 5^{i-j}$ for each $i\geq j$.
We conclude that $\tau(x')\neq \tau(x)$.
\qed

\subsection{More existence results of Lyapounov functions}

We will use the following easy Lemma in our next result:

\begin{lem}\label{product}
Let $\tau_i, 1\leq i \leq k,$ be finitely many  non negative Lyapounov functions,
then the product $\tau=\tau_1\tau_2 \cdots \tau_k$ is a non negative  Lyapounov function.
If all the $\tau_i$ are regular at some point $x_0$, then so is $\tau$.
\end{lem}

\proof
By  induction, it is enough to prove the statement for $k=2$.
The expression
$$
d\tau(x)=\tau_1(x) d\tau_2(x) +\tau_2(x)d\tau_1(x)
$$
implies that $d\tau_x\cdot v\geq 0$ for each $(x,v)\in \mC$.
Assume now that there exists $(x,v)\in \mC$, $v\neq 0$, such that $d\tau_x\cdot v=0$.
Then each of the terms $\tau_1(x) d\tau_2(x) \cdot v$ and $\tau_2(x)d\tau_1(x)\cdot v$ vanish,
which implies that each of the linear forms $\tau_1(x)d\tau_2(x)$ and $\tau_2(x)d\tau_1(x)$ vanish, hence that $d\tau(x)=0$.
We have proved that $\tau$ is a  Lyapounov function. If the functions $\tau_1$ and $\tau_2$ are regular at $x_0$,
then $\tau_i(x_0)>0$ and we see that $d\tau(x_0)\neq 0$.
\qed

For $K\subset M$, we set $\mF_{\mC}^{\pm}(K):=\bigcup _{x\in K}\mF_{\mC}^{\pm}(x) $.

\begin{prop}\label{proplyapnulcomp}
	Let $K\subset M$ be a compact set.
	Then there exists a  non negative Lyapounov function $\tau_+$ such that $\tau_+=0$
	on $K$ (hence on  $\mF^-_{\mC} (K)$) and $\tau_+>0$ outside of $\mF^-_{\mC} (K)$.
	This implies in particular that $\mF_{\mC}^-(K)$ is closed.
	The function $\tau_+$ can be chosen regular on 
	$\mD(\mC)-\big(\mF^-_{\mC}(K)\cup \mR_{\mC} \big)
	$.
	
	There also exists a non positive  Lyapounov function $\tau_-$ 
	such that $\tau_-=0$
	on $K$ (hence on  $\mF^+_{\mC} (K)$) and $\tau_-<0$ outside of $\mF^+_{\mC} (K)$.
 This implies that $\mF_{\mC}^+(K)$ is closed.
	The function $\tau_-$ can be chosen regular on 
	$\mD(\mC)-\big(\mF^+_{\mC}(K)\cup \mR_{\mC} \big)
	$.
\end{prop}

\proof
The second part of the statement is a consequence of the first part applied to the reversed cone $-\mC$.
More precisely, we have $\tau_-(\mC)=-\tau_+(-\mC)$.

To prove the first part, we  fix a point $x_0\in M-\mF^-_{\mC}(K)$.
For each $y\in K$, there exists a  Lyapounov function $f$ such that 
$f(y)<f(x_0)$. If moreover $x_0\not \in \mR_{\mC}$, then the function $f$ can be chosen regular at $x_0$.
By composing $f$ on the left with a non decreasing function,
we deduce the existence of a Lyapounov function $\tau_y$ such that $\tau_y \geq 0$, $\tau_y=0$ in a neighborhood $U_y$ of $y$, and $\tau_y(x_0)>0$.
If $x_0\not \in \mR_{\mC}$, then in addition $\tau_y$ is regular at $x_0$.

Since $K$ is compact, there exist finitely many points $y_1, \ldots, y_k$ such that the corresponding open sets $U_{y_i}$ cover $K$.
The product $\tau_0:=\tau_{y_1} \tau_{y_2} \cdots \tau_{y_k}$ is a  non negative Lyapounov function such that 
$\tau_0(x_0)>0$, and, if $x_0\not \in \mR_{\mC}$, $d\tau_0(x_0)\neq 0$.

For each $x_0\in M-\mF^-_{\mC}(K)$, we have proved the existence of an open neighborhood $V_0$ of $x_0$ and of a  non negative Lyapounov function $\tau$
which is null on $K$ and positive on $V_0$.
We can cover the separable metric space $M-\mF^-_{\mC}(K)$ by a sequence $V_i$ of open sets such that, for each $i$,
there exists a non negative Lyapounov function $\tau_i$ which is null on $K$ and  positive on $V_i$.
Then there exists a positive sequence $a_i$ such that $\tau:= \sum _i a_i\tau_i$ is a smooth non negative function which is positive on 
$M-\mF^-_{\mC}(K)$.

By exactly the same method we can also obtain a  non negative Lyapounov function $\tau$ which is null on $K$ and which has the property that 
$d\tau_x\cdot v>0$ for each 
$ x\in M-(\mF^-_{\mC}(K)\cup \mR_{\mC})$ and $v\in \mC(x)$.
\qed

By adding the functions $\tau_+$ and $\tau_-$, we obtain:

\begin{cor}\label{corlyapnul}
	Given a compact $K\subset M$, there exists a Lyapounov function 
	which is null on $K$ and regular on
	$\mD(\mC)-\big(\mF_{\mC}(K,K) \cup \mR_{\mC}\big)$, where 
	$\mF_{\mC}(K,K):=\mF_{\mC}^+(K)\cap \mF_{\mC}^-(K)$.
\end{cor}

Let us also state the following :

\begin{prop}\label{proplyapnul}
	Let $A\subset M$ be a trapping domain.
	There exists a  Lyapounov function $\tau$ such that $\tau>0$ on $A$ and $\tau<0$ outside of $\bar A$.
	The function $\tau$ can be chosen regular on 
	$\mD(\mC)- (\mR_{\mC}\cup \partial A) 
	$.
\end{prop}

\proof
We consider an enlargement $\mE$ of $\mC$ such that $A$ is 
a trapping domain  for $\hat \mE$.

We first fix a point $x_0\in A$ and prove the existence of a  Lyapounov
function which is non negative, null outside of $A$, positive at $x_0$ and, if $x_0$ is not stably recurrent, 
regular at $x_0$.

We consider a point $x_1\in A\cap \mI^-_{\mE}(x_0)$.
Then the set $A_1:=\mI^+_{\mE}(x_1)$ is open, it contains $F_i:=\mF^+_{\mC}(x_0)$,
and its closure is contained in $\mF^+_{\hat \mE}(x_1)$, hence in $A$.
In other words, the closure of  $A_1$ is disjoint from the {set} $F_e:=M-A$.
By Proposition \ref{prop-Areg}, there exists a smooth (near $\mD(\mC)$) trapping domain 
$A'_1$ which contains $F_i$ and whose closure is disjoint from $F_e$.
By Proposition \ref{prop-fol}, there exists a  Lyapounov function $\tau: M\lto [-1,1]$ (for $\mC$)
which is equal to $1$ on $F_i$ and to $-1$ on $F_e$. The non negative Lyapounov function $1+\tau$
is then null outside of $A$ and positive at $x_0$.

In the case where $x_0$ is not stably recurrent and non degenerate, we can take $\mE$ in such a way that $x_0\not \in A_2:=\mI^+ _{\mE}(x_0)$,
hence $x_0$ belongs to the boundary of this trapping domain. The closure of $A_2$ is disjoint from the complement $F_e$
of $A$.
By Propositions \ref{prop-Areg} and \ref{prop-fol}, we find a non negative Lyapounov function $\tau$ which is regular (hence positive) at 
$x_0$ and null outside of $A$.

By considering a convex combination of countably many of the Lyapounov functions
we just built, we obtain a non negative Lyapounov function $\tau_i$ which
is positive on $A$ and regular on $(A \cap \mD(\mC))-\mR_{\mC}$.

We can apply the same result to the cone $-\mC$ and the trapping domain 
$M-\bar A$, and get a Lyapounov function $\tau_e$ (for $\mC$) which is non positive,
negative outside of $\bar A$, and regular on $(\mD(\mC)-\bar A)-\mR_{\mC}$.

The sum $\tau:= \tau_i+\tau_e$ then satisfies the conclusions of the 
proposition.
\qed

\subsection{Hyperbolic cone fields}\label{sec-hyp-proof}

We prove Theorem \ref{thm3} and discuss some alternative characterizations of global hyperbolicity.
We start with an easy observation:

\begin{lem}\label{lemclosed}
If the closed cone field $\mC$ 
satisfies \ref{GH2}, then $\mathcal{J}_{\mC}^\pm (x)$ is closed for each 
$x\in M$. 
\end{lem}

{The lemma shows that  hyperbolic cone fields satisfy the analogous conditions to 
causal simplicity in \cite{ms1}.}

\proof
  Let $y_n\in \mathcal{J}_{\mC}^+(x)$ be a convergent sequence 
with limit $y\in M$. 
Let $Y$ be the compact set $Y:= \{y,y_1,y_2, \ldots\}$.
The set $\mJ_{\mC}(x,Y)$ is compact and it contains $y_n$ for each $n$,
hence it contains the limit $y$.
\qed

Let us denote by $\mC_K$ the cone field which is equal to $\mC$ on $K$ 
and degenerate outside of $K$.
If $\mC$ is a closed cone field and $K$ is a  closed set, then $\mC_K$ is a closed cone field. If $\mC$ is causal, then so is $\mC_K$. 

\begin{lem}\label{lemL}
Let $\mC$ be a causal closed cone field and $K$ be a compact set.
Then there exists an open enlargement $\mE$ of $\mC_K$ and a real number
$L>0$ such that all $\mE$-timelike curves have length less than $L$.
	\end{lem}

\proof
Let $\mE_n$ be a decreasing sequence of open cone fields converging to $\mC_K$.
We can assume that $U_n:=\mD(\mE_n)$ is bounded for each $n$.
 If the conclusion of the Lemma does not hold, there  exists a sequence 
 $\gamma_n:[-l_n,l_n]\lto M$ of $\mE_n$-timelike curves parametrized by arclength
 with $l_n$ unbounded.
 By Proposition \ref{proplim}, there  exists a complete $\mC_K$-causal
curve $\gamma:\Rm\lto M$. Since $\mC_K$ has no singular points, this
curve has infinite length in the forward direction.
  Let $\omega$ be a limit point of $\gamma$ at $+\infty$.
  For each $s>t\in \Rm$, we have $\gamma(s)\in \mJ^+_{\mC}(\gamma(t))$. Since this set is closed (Lemma \ref{lemclosed}),
  we deduce that $\omega \in \mJ^+_{\mC}(\gamma(t))$, or in other words that $\gamma(t)\in \mJ^-_{\mC}(\omega)$, and this holds for all $t$.
  Since $\omega$ is not singular, there exists a local time function,
  and this implies that $\gamma$ has another limit point $\omega'$.
  Since $\mJ^-_{\mC}(\omega)$ is closed, we obtain that $\omega'\in  \mJ^-_{\mC}(\omega)$, and 
  similarly  $\omega\in  \mJ^-_{\mC}(\omega')$.
  This is in contradiction with $\mC$
  being causal.
  \qed
 
\begin{cor}
Let $\mC$ be a  hyperbolic closed cone field and $K$ be a compact set.
	The stably recurrent set $\mR(\mC_K)$ is empty. 
\end{cor}

\proof
If $\mR(\mC_K)$ is not empty, then $\mC_K$ has a complete causal curve, by Corollary \ref{cor-R}.
Since $\mC_K$ has no singular point, this curve has infinite length, 
which contradicts Lemma \ref{lemL}.
\qed

\begin{cor}\label{cor52}
	Let $\mC$ be a hyperbolic closed cone field, and $K_1,K_2$
	be two compact sets. Let $K$ be a compact set containing 
	$\mJ_{\mC}(K_1,K_2)$. Then
	$$
	\mF_{\mC_K}(K_1,K_2)=\mJ_{\mC_K}(K_1,K_2)=\mJ_{\mC}(K_1,K_2).
	$$
\end{cor}

\proof
The second equality is clear.
To prove the first equality, we consider a  sequence
$\mE_n$ of open enlargements of $\mC_K$ decreasing  to  $\mC_K$.
By Lemma \ref{lemL}, we can assume that each $\mE_1$-timelike curve has length
less than $L>0$. This is then true for all $\mE_n$.
Given $x\in \mF_{\mC_K}(K_1,K_2)$, there exists a sequence 
 $\gamma_n:[0,1]\lto M$ 
of $\mE_n$-timelike curves connecting $K_1$ to $K_2$,  parametrized 
proportionally to arclength, and passing through $x$. Since the curves $\gamma_n$ have bounded
length, they are equi-Lipschitz. Up to a subsequence, they converge 
uniformly to
a Lipschitz curve $\gamma:[0,1]\lto M$ which is $\mC_K$-causal by Lemma 
\ref{lemlength},  passes through $x$, and connects $K_1$ to $K_2$. 
This implies that $x\in \mJ_{\mC_K}(K_1,K_2)$.
\qed

\textsc{Proof of Theorem \ref{thm3}:}
We first prove the existence of a steep Lyapounov function for a  hyperbolic cone 
field $\mC$.
Let $K_i, i\geq 0$ be a sequence of compact subsets of $M$  such that $\mJ_{\mC}(K_i,K_i)$
is contained in the interior of $K_{i+1}$ and such that $M=\cup_i K_i$.
We set $A_i:= K_i\cap \mD(\mC)$.

For each $i \geq 2$, we apply 
Corollary  \ref{corlyapnul} 
to the cone field $\mC_{K_{i}}$ and the compact set $K_{i-2}$.
Since $\mF_{\mC_{K_{i}}}(K_{i-2},K_{i-2})=\mJ_{\mC}(K_{i-2},K_{i-2})\subset \mathring K_{i-1}$ by Corollary \ref{cor52}
and since $\mR(\mC_{K_{i}})$ is empty, we obtain 
 a smooth function $\tau_i:M\lto \Rm$ 
with the following properties:
\begin{itemize}
	\item[$\cdot$] $\tau_i$ is a Lyapounov function on $K_{i}$, which means that 
	$d\tau_i(x)\cdot v>0$ for each $x\in K_{i}$ such that $d\tau_i(x) \neq 0$
	and each $v\in \mC(x)$.
	\item[$\cdot$] $\tau_i$ is regular on $A_{i}-\mathring K_{i-1}$, which means that $d\tau_i(x)\neq 0$ for each $x\in 	A_{i}-\mathring K_{i-1}$.
	\item[$\cdot$] $\tau_i$ is null on $K_{i-2}$.
\end{itemize}
We also let $\tau_1$ be a smooth function on $M$ which is a Lyapounov function on $K_1$ and regular on $A_1$.

We now prove the existence of  a sequence $a_i$ of positive numbers such that 
the sum $\tau:= \sum_{i\geq 1} a_i \tau_i$ is a steep Lyapounov function.
Note that this sum is locally finite.

We build the sequence $a_i$ by induction, in such a way that the partial sum
$\sum_{i=1}^k a_i \tau_i$ is a steep Lyapounov function on $K_{k}$ for each $k$.

The function $\tau_1$ is a Lyapounvov function on the compact set 
$K_1$, hence there exists $a_1>0$ such that $a_1\tau_1$ is steep on 
$K_1$.
The function $\tau_2$ is  Lyapounov on $K_2$ and regular on 
$A_2-\mathring K_1$. Then there exists $a_2>0$ such that $a_1\tau_1+a_2\tau_2$
is a steep Lyapounvov function on $A_2-\mathring K_1$, hence on $K_2$
(being steep is an empty condition outside of $\mD(\mC)$).
Assuming that $a_1, \ldots, a_k$ have been constructed, observe that the
function $\tau_{k+1}$ is Lyapounov on $K_{k}$ and regular
on $A_{k}-\mathring K_{k-1}$. On the other hand the partial sum 
$\sum_{i=1}^k a_i\tau_i$ is 
a smooth function on $M$ which is a steep Lyapounov function on $K_{k}$.
There exists $a_{k+1}>0$ such that $\sum_{i=0}^{k+1} a_i\tau_i$
is a steep Lyapounvov function on $A_{k+1}-\mathring K_{k}$, hence on $K_{k+1}$.
This ends the proof of the existence of a steep Lyapounov function.

Conversely, let us assume the existence of a steep Lyapounvov function $\tau$.
It is clear that $\mC$ is causal.
Let us prove that $\mF^{\pm}_{\mC}(x)=\mJ^{\pm}_{\mC}(x)$ for each $x$. 
We consider a decreasing sequence $\mE_n$ of enlargements of $\mC$, which have the property that $d\tau_y\cdot v \geq |v|_y/2$ for each $(y,v)\in \mE_n$.
Given $z\in \mF^+_{\mC}(x)$, there exists a sequence $\gamma_n:[0,1]\lto M$ of
smooth  $\mE_n$-timelike
curves such that $\gamma_n(0)=x$ and  $\gamma_n(1)=z$.
We can assume that $\gamma_n$ is parametrized proportionally
to arclength, hence is $L_n$-Lipschitz, where $L_n$ is the length of $\gamma_n$.
The hypothesis made on $\mE_n$ implies that $L_n\leq 2(\tau(z)-\tau(x))$
is bounded.
At the limit, we obtain a Lipschitz causal curve 
$\gamma:[0,1]\lto M$ connecting $x$ to $z$. We have proved that 
$\mF^+_{\mC}(x)\subset \mJ^+_{\mC}(x)$, hence these sets are equal.

We finally prove  \ref{GH2}.
The set $\mJ_{\mC}(K,K')=\mF_{\mC}(K,K')$ is closed.
If $\gamma$ is a causal curve joining $K$ to $K'$, then the length of $\gamma$
is bounded by $\max_{K'}\tau-\min _K \tau$.
This means that $\gamma$ is contained in a bounded set, hence that 
$\mJ_{\mC}(K,K')$ is bounded.
Being closed and bounded in the complete Riemannian manifold $M$,
the set $\mJ_{\mC}(K,K')$ is compact.
\qed

Let us finish this section with some alternative characterizations of global hyperbolicity
which generalizes \cite[Theorem 3.79]{ms1} to the present case:

\begin{prop}\label{prop-prop}
	A closed cone field $(M,\mC)$ is  hyperbolic if and only if it is  regular, and if
	\begin{itemize}
		\item[\mylabel{GH4}{(GH4)}] 
		For each compact $K\subset M$, there exists $L>0$ such that each causal curve $\gamma:[0,T]\lto M$ 
		of length more than $L$ satisfying $\gamma(0)\in K$
		satisfies $\gamma(T)\not \in K$.
	\end{itemize}
\end{prop}

\proof
By Theorem \ref{thm3} there exists a steep Lyapounov function $\tau$ if $(M,\mC)$ is hyperbolic.
This implies 
\ref{GH4} with $L= \max_K \tau-\min_K \tau$.

Conversely, assume that $\mC$ is a  regular closed cone field satisfying  \ref{GH4}.
Then $\mC$ is causal (it satisfies \ref{GH1}). 
Moreover, if $K,K'\subset M$ are  compact, there exists an upper bound on 
the length of causal curves with endpoints in $K\cup K'$, by \ref{GH4}.  In view of Lemma \ref{lem-limit},
this implies  \ref{GH2}.
\qed

Let us discuss the case where $V$ is a smooth vector field without singular points generating a complete flow $\phi^t(x)$.
We say that the dynamics is trivial if there exists a submanifold transverse to $V$ and intersecting
each orbit in one and only one point. This is equivalent to the existence of a steep Lyapounov function
for $V$.

We say that the action of $\phi$ is proper if, for each compact $K\subset M$, there exists $L>0$ such that $\phi^t(x)\not \in K$
if $x\in K$ and $t>L$. This is equivalent to \ref{GH4}.
Proposition \ref{prop-prop} implies
that the dynamics is trivial if and only if the action  $\phi$ of $\Rm$ on $M$ is  proper.

We now  give  another characterization of global hyperbolicity in the spirit of the one given by
	Minguzzi for Lorentzian metrics, \cite{mi09}.
	
	\begin{prop}
		A closed cone field $(M,\mC)$ is  hyperbolic if and only if 
		\begin{itemize}
			\item[\mylabel{GH5}{(GH5)}] No complete causal curve is contained in a compact set.
			\item[\mylabel{GH6}{(GH6)}] for all $K,K'\subset M$ compact the set $\mJ_{\mC}(K,K')$ is bounded.
		\end{itemize}
	\end{prop}

	\proof
	It is easy to see that the existence of a steep Lyapounov function implies \ref{GH5} and \ref{GH6}.

	Conversely assume that the closed cone field $(M,\mC)$ satisfies  \ref{GH5} and \ref{GH6}. \ref{GH5} implies that there are no singular points.
	We now prove \ref{GH4}.
	Consider a compact subset $K$, and assume that \ref{GH4} does not hold:
	There exists a sequence $\gamma_n\colon [-b_n,b_n]\lto M$ of causal curves with boundaries in $K$, parametrized by arclength, with  $b_n\lto \infty$.
	Then by \ref{GH6} $\gamma_n$ is contained in the bounded set $\mJ_{\mC}(K,K)$.
 By Proposition \ref{proplim} a subsequence converges 
	to a  complete causal curve $\gamma\colon  \Rm\lto M$ which is contained in the compact  set  $\overline{\mJ_{\mC}(K,K)}$.
	This contradicts \ref{GH5}.
	\qed

\section{Final remarks on  the stably recurrent set}

We propose here some additional remarks on the stably recurrent set $\mR_\mC$.
We first improve Corollary \ref{cor-R}:

\begin{prop}
	Let $\mC$ be a closed cone field, and $\mR_{\mC}$ be the stably recurrent set.
	For each $x\in\mR_\mC$, there exists a complete causal curve $\gamma$ which takes values in $\mR_\mC$ and satisfies 
	$\gamma(0)=x$.
\end{prop}

\proof
 Let $\tau_0$ be a Lyapounov function which is regular outside of $\mR_{\mC}$.
We consider a decreasing sequence $\mE_n$ of open enlargements of $\mC$ all smaller than $\{d\tau_0>0\}$.
As in Corollary  \ref{cor-R},  let $\gamma_n:\Rm\lto M$ be a sequence  of $\mE_n$-timelike periodic curves 
parametrized by arclength and satisfying $\gamma_n(0)=x$.
The function $\tau_0$ is non decreasing, hence  constant, on $\gamma_n$.
At the limit, we obtain a complete causal curve $\gamma$, and the function $\tau_0$ is constant on it.
We deduce that the curve $\gamma$ takes values in the critical set of $\tau_0$, i.e. in $\mR_\mC$. 
\qed

We finish with a stability property:

\begin{prop}\label{semicontinu}
		Let $\mC$ be a closed cone field. We assume that  the stably recurrent set $\mR_{\mC}$ is compact.
	Then for every neighborhood $U$ of $R_\mathcal{C}$ there exists
	a closed enlargement $\mC_U$ of $\mC$ such that 
	$\mR_{\mathcal{C}_U}
	\subset U$.
\end{prop}

\proof
We assume, without loss of generality, that $U$ is bounded, hence $\partial U$ is compact.
It is enough to prove the existence of an enlargement $\mE$ of $\mC$ such that 
$\mR_{\hat \mE}$ is disjoint from $\partial U$.

Let us fix a point $z\in \partial U$. By Lemma  \ref{recreg}, there exists a Lyapounov
function $\tau^z$ for $\mC$ such that $a:= \tau^z(z)$ is a regular value of 
$\tau^z$. Then, there exists a closed enlargement $\mC^z$ of $\mC$ such that 
$\tau^z$ is a regular Lyapounov function for $\mC^z$ in a neighborhood
of $\{\tau^z=a\}$. This implies, by Lemma \ref{lem-trap}, that $\{\tau^z>a\}$
is a trapping region for $\mC^z$, hence that $z\not \in \mR_{\mC^z}$.

The open sets $M-\mR_{\mC^z}$, $z\in \partial U$, thus cover the compact set $\partial U$,
hence finitely many of them cover $\partial U$.
By taking the intersection of the corresponding cone fields $\mC^z$, 
we obtain a closed enlargement of $\mC$ whose stably recurrent set is disjoint from 
$\partial U$, as was claimed.
\qed

\end{document}